\documentclass[11pt]{article}
\usepackage[english]{babel}
\usepackage[latin1]{inputenc}
\usepackage{amsfonts,amssymb,amsmath, epsfig}
\usepackage{color,graphicx,graphics}
\usepackage{amsmath,amstext,amssymb,amsfonts, amscd}
\usepackage{hyperref}
\usepackage{a4wide}
 \numberwithin{equation}{section}
\def\Box{\leavevmode\vbox{\hrule
     \hbox{\vrule\kern4pt\vbox{\kern4pt}%
           \vrule}\hrule}}
\def\blackbox{\leavevmode\vrule height 5pt width 4pt depth 0pt\relax}
\def\endproof{\null\hfill {$\blackbox$}\bigskip}

\newcounter{appendix}
\def\appendix{\advance\c@appendix by 1
   \def\thesection{\Alph{section}}
   \ifnum\c@appendix=1 \setcounter{section}{-1} \fi
   \@startsection {section}{1}{\z@}{-3.5ex plus -1ex minus
   -.2ex}{2.3ex plus .2ex}{\Large\bf}}

\newtheorem{lemma}{Lemma}[section]

\newtheorem{proposition}[lemma]{Proposition}

\newtheorem{remark}{Remark}[section]

\begin {document}
\title{Kinetic Theories for Metropolis Monte Carlo Methods  }
\author{Michael Herty \footnote{RWTH Aachen University, Institute Geometry and Practical Mathematics, 52056 Aachen, Germany} \and  
	Christian Ringhofer\footnote{Arizona State University, School of Mathematical and Statistical Sciences, Tempe, AZ 85287-1804, USA} }
\date{\today}
\maketitle
\begin{abstract}
We consider generalizations of the classical inverse problem to Bayesien type estimators, where the result is not one optimal parameter but an optimal probability distribution in parameter space.
The practical computational tool to compute these distributions is the Metropolis Monte Carlo algorithm.
	We derive kinetic theories for the Metropolis Monte Carlo method in different scaling regimes. The derived equations yield a different point of view on the classical algorithm. It further inspired modifications to exploit the difference scalings shown on an simulation example of the Lorenz system. 
\end{abstract}


\section {Introduction}\label {sintro}
The subject of this paper is the solution of  the following problem.
Given a model depending on a set of parameters and a set of observed data, find the optimal parameters to fit the given data.
In the framework of classical inverse problems this results into the minimization problem
\begin{equation} \label {lminx}
dist(\nu ,\mu (x))\rightarrow  \min _x \ ,
\end{equation}  
where  $\nu = \{ z_i \in \mathbb{R}^{N_{\nu}}: i \in \mathbb{N} \} $ denotes the set of observed data $z_i$ and $\mu(x): \mathbb{R}^{N_x} \to   \mathbb{R}^{N_{\nu}} $ is  obtained by the  model $\mu$ using the (often high but finite--dimensional) parameter $x \in \mathbb{R}^{N_x}$. The function
$dist (\nu ,\mu )$ denotes some measure of the distance between the measures $\nu $ and the model  $\mu $ (usually some norm of the form $dist(\nu ,\mu )=  \sum_{i=1}^{|\nu|}  || z_i -\mu ||$).
In practical applications, the problem (\ref {lminx}) often turns out to be extremely ill conditioned, exhibiting multiple optima, and has to be solved via various regularization techniques, see e.g. \cite{MR1009037,MR4619819,MR4565586}.


In this paper, we will consider the more general approach, seeking not a single optimal parameter value $x$ but a probability distribution $P(x)$ on the space of parameters \cite{Sullivan2015,Kaipio2005,Somersalo2007}.
So, given a set of observed data $\nu$, we compute a probability distribution $P^\nu $ of the observed data.
The goal is to compute a probability distribution $P(x)$ in the parameter space and a corresponding distribution $P^\mu [P]$ in the data space given by
\begin{align*}
P^\mu [P](z)=\int  \delta (z-\mu  (x))  dP(x) \ ,
\end{align*}
with $\delta-$being the Dirac-delta distribution, such that a distance
\begin{equation} \label {lminP}
dist (P^\nu ,P^\mu [P]) \to  \min _P
\end{equation}  
becomes minimal. Here, $dist(*)$ denotes a measure of the distance between the measures $P^\nu$ and $P^\mu.$ Further, $\mu (x)$ denotes still the model  computed  for a single parameter $x $ and $P^\mu [P]$ denotes the corresponding probability distribution of the results  in the case when  parameter $x$  is also distributed according to $P.$ 
We note that the optimization problem (\ref {lminP}) reduces to the deterministic inverse problem (\ref {lminx}) if we reduce  the probability distribution $P$ to a degenerate distribution (i.e.,  a $\delta -$ distribution)  and compatible distances. 
Therefore, the optimization problem for a deterministic parameter $x$ can be embedded into an optimization problem for a probability density $P$.
\par 
The solution of this problem \eqref{lminP} yields significantly more information.
The optimal parameter $x_*$ will be chosen as the expectation $\mathbb{E}_P(X)$ where $X$ is a random variable with distribution $P$.  One also obtains additional information about the reliability of this parameter by considering the variance $\mathbb{V}_P(X)$ or the relative reliability of individual components of $x$ by computing the covariance $COV_P(X)$ in the case of a higher dimensional parameters.  Multiple local extrema of the classical inverse problem (\ref {lminx}) may show up as local peaks in the distribution $P(x)$ in the solution of the problem (\ref {lminP}).  
\par 
To actually solve the optimization problem (\ref {lminP}) we require sampling. We consider a special class of Markov chain Monte Carlo (MCMC) methods. Namely, we will focus on one of the simplest methods, the  Metropolis Monte Carlo (MMC) algorithm \cite{MR3363437}. Certainly, there is a vast literature available on MCMC and MMC methods available and we refer to e.g.  \cite{brooks2011handbook,MR2858443,MR2884617,MR2260716,MR4196544} for further references. Our focus here, is to derive meanfield limits and corresponding macroscopic equations. We focus on the simplest form of MMC to illustrate the ideas:  The following  algorithm produces a sequence of probability distributions $P_n$ which (hopefully) converge to a density $P_\infty $ solving the problem (\ref {lminP}).
\par 
MMC algorithm  allows for iterative updates on the data set $\nu $.  In many applications, such as in geophysics \cite{MR3126997} or in meteorology \cite{MR4687001},  the data set $\nu = \{z(t_k): k=1,\dots \}$ are updated at the same time $t_k$ the model  parameters or their distribution  are computed in real time. 
This generalizes the minimization problem (\ref {lminP}) to the time dependent problem
\begin{equation} \label {lminPt}
dist (P^{\nu }(t_k) ,P^\mu [P(t_k)])\rightarrow  \min _{P(t_k)},\
k=1,2,... \ ,
\end{equation}  
which is easily  incorporated in the iterative MMC algorithm.
\par 
\noindent {\bf Relation to Bayesian estimation}:
The presented approach is closely related to the methodology of Bayesian inversion. This relies  on  Bayes' formula
\begin{equation} \label {lBay}
\mathcal{P} (x|\nu )=\frac {\mathcal{P} (\nu |x)\mathcal{P} _{pri}(x)}{\int  \mathcal{P} (\nu |z)\mathcal{P} _{pri}(z)\ dz} \ .
\end{equation}   
Here, $\mathcal{P} (x|\nu )$ is the conditional probability of the parameter $x$, given the data distribution $\nu $ (the resulting distribution of the parameter $x$).
$\mathcal{P} (\nu |x)$ denotes the conditional probability of the actual data given the parameter $x$, which is modeled by $P_\mu $ in formulation (\ref {lminP}).
Finally, $\mathcal{P} _{pri}$ is the prior  distribution of the  parameter $x$ (the prior).
A `good' prior turns out to be all important for the success of Bayesian estimation, see e.g. \cite{MR3126997,MR2652785,MR4604099}. In practice, the resulting probability distribution $\mathcal{P} (x|\nu )$ in Bayes' formula (\ref {lBay}) may also be computed by the MMC.
So, the methodology in this paper can be interpreted as an iterative application of Bayesian estimation,
where the prior $\mathcal{P} _{pri}(x)$ is updated using the last iterate of the result $\mathcal{P} (x|\nu )$.
\par
This paper is organized as follows.
In Section \ref {sMMC} we define the general iterative Metropolis Monte Carlo algorithm studied in this paper, giving in iteration $n$ a distribution function $f_n(x)$ in the parameter space.
We derive a kinetic equation for the parameter distribution $f_n(x)$.

Section \ref {scntlim} is devoted to the behavior of the functions $f_n$ for many iterations, i.e. the convergence properties of the iterative MMC method.
Here, we will consider two regimes:\
In the first regime we consider the case of a small acceptance probability in the MMC algorithm but an arbitrary prediction probability. That is, we take a more or less arbitrary guess but accept the guess only infrequently. This leads to a Boltzmann type equation for the parameter density whose convergence properies wil be investigated using entropy methods, standard in kinetic theory.
In the second regime we use only predictions which only are a small random variation from the previously computed  data point, but use a more or less arbitrary acceptance probability in the MMC algorithm.
This leads to a Brownian motion regime and a  Fokker - Planck type equation for the parameter distribution.    

In Section \ref {sec:micromacro} we use the results of Section \ref {scntlim} to accelerate the convergence of the MMC algorithm by decomposing the parameter density $f_n$ into a part which is computed by the classical MMC algorithm and a part where we only use the macroscopic models from Section \ref {scntlim} (which are cheaper to evaluate). Finally, in Section \ref {slorentznum} we verify the above results on a simple model with a three dimensional parameter space.
We choose the  chaotic Lorentz system taken from \cite{MR3363508}.

The more technical proofs and the details of the numerical implementation 
are deferred to the Appendix in Section \ref {sappnd}.

\section {The Metropolis Monte Carlo algorithm (MMC)}\label {sMMC}
In this section we will define the general MMC algorithm and analyze  its convergence properties.
We first take the point of view that we want to compute a given distribution $P(x)$ by computing a sequence of particles $x_n,\ n=1,2,...$.

\subsection {General definition of the method}\label {ssgendef}
The general structure  of the MMC  algorithm is the following: 
\\ \\ 
Given a probability distribution $P_n(x)$ computed from the particles $x_1,..,x_n$:
\vspace{1cm}
\\{\bf  Step 1:
Compute a random proposal $x_p $ from a probability distribution $\tau $.
\vspace{5mm}
\\ Step 2: 
Compute an acceptance rate $\alpha \in [0,1]$ and accept the proposal $x_p$ with a probability $\alpha $ or reject the proposal $x_p $ with probability $\mathcal{P} =1-\alpha $.
\vspace{5mm}
\\Step 3: Compute $P_{n+1}$ by either adding  $x_{n+1}=x_p$ or $x_{n+1}=x_n$ to the distribution. 
}
\vspace{1cm}
\\ 
The method will either add the proposal $x_p$ to the distribution or, in the case of rejection, reinforce the parameter $x_n$ by adding an additional copy to the distribution $P_n.$ 
The proposal particle $x_p$ and the acceptance rate $\alpha $ will be chosen dependent on the already computed values $x_1,..,x_n$ and the given distribution of the observed data distribution $P^\nu $ at iteration $n.$  In particular, the acceptance rate $\alpha $ will be chosen dependent on whether adding $x_p$ or another copy of $x_n$ to the distribution will make the modeled data distribution $P^\mu  [P_{n+1}]$ a better match to the observed data distribution  $P^\nu $ in the sense of the distance $dist(*).$  In this paper we will make the proposal $x_p$ and the acceptance rate $\alpha $ dependent not only on the last computed node $x_n$ but also on the whole up to date computed distribution $P_n$.
\par 
We will restrict the dependence on $P_n$ to the dependence on a certain number of moments of $P_n$, i.e. considering quantities like means or variance.  We define the moments  $\kappa  _n \in \mathbb{R}^{N_\kappa}$ of the distribution $P_n$ as
\begin{equation}\label{kappa} 
\kappa  _n=\frac 1 n \sum _{j=1}^n (x_j,x_j^2,..)  \ .
\end{equation}
So, the general MMC algorithm considered in the following is of the form
\begin{equation} \label {genmmc}
\hbox {General MMC algorithm}
\end{equation} 
Given the discrete distribution $P_n$by the already computed nodes $x_1,..,x_n$ and the moments of $P_n$ by $\kappa  _n=\frac 1 {n}\sum _{j=1}^{n}(x_j,x_j^2,..)$:
\begin {itemize}
\item [{\bf Step 1}]  (Proposal) Compute a random proposal $x_p $ from the probability distribution $\tau $:
$$
d\mathcal{P} [x_p =z]=\tau (z|x_n,\kappa  _n)\ dz 
$$
\item [{\bf Step 2}] (Acceptance/rejection step)
To compute the acceptance probability $\alpha $ we have to compare the quality of $P_{n+1}$ if we either add $x_p$ or another copy of $x_n$ to the distribution.
The update of the moment vector $\kappa  _n$ is of the form
$$
\kappa  _{n+1}(x,\kappa  _n)=\frac {\Delta  \kappa  (x)+n\kappa  _n}{n+1},\
\Delta  \kappa  (x)=(x,x^2,..)
$$
with either $x=x_p$ or $x=x_n$.
So, we compute
$\alpha (x_p,\omega _n(x_p,\kappa  _n),x_n,\omega _n(x_n,\kappa  _n))$ with
$$
\omega _n(x,\kappa  )=\frac {\Delta  \kappa  (x)+n\kappa  }{n+1},\
\Delta  \kappa  (x)=(x,x^2,..)
$$
\item [{\bf Step 3}](Update)
Set 
$$
x_{n+1}=x_p,\ \kappa  _{n+1}=\omega _n(x_p,\kappa  _n) \ \hbox {with probability} \ 
P= \alpha (x_p,\omega _n(x_p,\kappa  _n),x_n,\omega _n(x_n,\kappa  _n))
$$
or, set
$$
x_{n+1}=x_n,\ \kappa  _{n+1}=\omega _n(x_n,\kappa  _n) \ \hbox {with probability} \ 
P= 1-\alpha (x_p,\omega _n(x_p,\kappa  _n),x_n,\omega _n(x_n,\kappa  _n))
$$
\item [{\bf Step 4}] (Repeat) $n\rightarrow n+1$
\end {itemize}
So $x_{n+1}$ equals the proposal $x_p $ if $x_p $ is `better' than $x_n$ (determined by the acceptance rate $\alpha $)
or $x_{n+1}=x_n$ otherwise.
{\remark \label {rdep}
The acceptance rate $\alpha $ (and possibly the proposal probability $\tau$) may depend on the observed data distribution $P^\nu $. Since $P^\nu $ is taken as a given within the MMC algorithm we suppress this dependence in the notation for convenience.
}
{\remark\label {rtimdep}
There are two distinct scenarios for the estimation of the optimal parameter distribution.

In the first scenario all data $\nu $ and the distribution $P^\nu $ are known during the execution of the algorithm  (\ref {genmmc}).

In the other scenario, the data $\nu $ are updated continuously while the algorithm  (\ref {genmmc}) is executed. 
This scenario applies to c.f. weather prediction where the model parameters and the resulting predictions are updated as new data arrive.
\cite {weather1}.

In this scenario the acceptance rates $\alpha $ would have to be updated within the algorithm (\ref {genmmc}). Due the implicit dependence of the acceptance rate on the observed data distribution $P^\nu $ this would result in $\alpha $ being dependent on the iteration index of the MMC algorithm ,
so $\alpha =\alpha _n(x_p,\omega (x_p,\kappa  _n),x_n,\omega (x_n,\kappa  _n))$ in Step 3 of algorithm (\ref {genmmc}) would hold.
}

{\remark \label {rss}
In general the start up phase of algorithm (\ref {genmmc}), say the first 1000 nodes $x_n$, will be discarded to compute a steady state in algorithm (\ref {genmmc}), and using only  $N-1000$ nodes overall for the distribution $P(x)$~ \cite {startup1,startup2}.
}

 \subsection {Evolution of the probability density}\label {ssdensev}
 
 In this section , we derive the evolution equation for the probability density
$$ 
f_n(x',\kappa ' )\ dx'\kappa ' =d\mathcal{P} [x_n=x',\kappa _n=\kappa '] \ .
$$
corresponding to the Metropolis Monte Carlo algorithm (\ref {genmmc}).

For the derivation of the limiting density $f_\infty (x,\kappa )=\lim _{n\rightarrow \infty}f_n(x,\kappa )$ it will be convenient to derive the evolution equation for $f_n$ in  weak form.
Summing over all possibilities in algorithm (\ref {genmmc}), we have
\begin{align*}
&  f_{n+1}(x,\kappa ) = \\
&   \int  \delta (x_p-x)\d (\omega _n(x_p,r)-\kappa )\alpha (x_p,\omega _n(x_p,r),y,\omega _n(y,r))\tau (x_p|y,r)f_n(y,r) \ dx_pyr+ \\
  & \int  \delta (y-x)\d (\omega _n(y,r)-\kappa )[1-\alpha (x_p,\omega _n(x_p,r),y,\omega _n(y,r))]\tau (x_p|y,r)f_n(y,r)\ dx_pyr \ ,
  \end{align*}
or, in weak form,
\begin{align*}
&	\int  \phi (x,\kappa )f_{n+1}(x,\kappa )\ dx\kappa =\\
	& \int  \phi (x,\omega _n(x,r))\alpha (x,\omega _n(x,r),y,\omega (y,r))\tau (x|y,r)f_n(y,r) \ dxyr +\\
	& \int  \phi (y,\omega _n(y,r)) [1-\alpha (x_p,\omega _n(x_p,r),y,\omega _n(y,r))]\tau (x_p|y,r)f_n(y,r)\ dx_pyr  \ ,
\end{align*}
for all test functions $\phi (x,\kappa ) \in C^\infty_0( \mathbb{R}^{N_x} \times \mathbb{R}^{N_\kappa}; \mathbb{R}  ).$
After renaming $x_p \rightarrow x$ in the second integral,
\begin{align}\label {lver3}
&\int  \phi (x,\kappa )f_{n+1}(x,\kappa )\ dx\kappa = \\
&\int  \phi (x,\omega _n(x,r))\alpha (x,\omega _n(x,r),y,\omega (y,r))\tau (x|y,r)f_n(y,r) \ dxyr + \\
& \int  \phi (y,\omega _n(y,r)) [1-\alpha (x,\omega _n(x,r),y,\omega _n(y,r))]\tau (x|y,r)f_n(y,r)\ dxyr. \label{lverf3b}
\end{align}
The evolution equation in its weak form (\ref {lver3}--\ref{lverf3b}) will be used in the following to derive the evolution of the probability density in various limits, where the index $n$ is replaced by a continuous variable. 
 
 
\section {Continuum limits for a large number of iterations}\label {scntlim}
In this section we will derive the continuum limit for a continuous index (i.e. for a large number $N$ of iterations) 
in two different regimes. We will essentially replace the discrete index $n$ in (\ref {lver3}) by a continuous limit \\ $s=nh,\ n=1:N,\ h=N^{c-1} ,\ s\in (0,N^c)$ for some constant $c$ with $0<c<1$. 
So, the limit $N\rightarrow \infty ,\ h\rightarrow 0$, with $N$ the total number of iterations performed, will correspond to a continuous variable $s\in [0,\infty )$.  
We define
\begin{align}
f(x,\kappa ,s)=f_n(x,\kappa ),\
s=nh,\ n=1:N,\ h=N^{c-1} ,\ s\in (0,N^c)\ .
\end{align}
and
\begin{align}
& \omega _n(x,\kappa )=\frac {\Delta \kappa (x)+n\kappa }{1+n}\Rightarrow  \omega (x,\kappa ,s)=\frac {h\Delta \kappa (x)+s\kappa }{h+s}=\kappa +h\omega _1(x,\kappa ,s),  \\
&
\label {lom1}
\omega _1(x,\kappa ,s) =\frac {\Delta \kappa (x)-\kappa }{h+s} \ .
\end{align}
So, the term $\omega _1(x,\kappa ,s)$ gives the incremental change in the moments from one iteration step to the next.
Thus, (\ref {lver3}) yields 
\begin {align}\label {lver4}
& \int  \phi (x,\kappa )f(x,\kappa ,s+h)\ dx\kappa = \\
& \int  \phi (x,r+h\omega _1(x,r,s))\alpha (x,r+h\omega _1(x,r,s),y,r+h\omega _1(y,r,s))\tau (x|y,r)f(y,r,s) \ dxyr +
\\ 
& \int  \phi (y,r+h\omega _1(y,r,s)) [1-\alpha (x,r+h\omega _1 (x,r,s),y,r+h\omega _1(y,r,s))]\tau (x|y,r)f(y,r,s)\ dxyr  \ .
\end{align}
{�\remark {\label {rmkc}
We note that, at first glance, the choice of the constant $c\in (0,1)$ seems to be   arbitrary and only serves to yield a continuous limit  problem for the continuous variable $s\in (0,\infty )$, and is just an interpretation of the discrete problem  (\ref {lver3}).
However, to obtain a continuum   limit for a total of $N$ iterations, we have to have that $h=N^{c-1}$ tends to zero.
Also, if we take only small steps  to compute the next proposed node in the distribution, as is the case in the Brownian motion regime in Section \ref {ssbrown}, we have to have that $N\mathbb E (\tau )=Nh=N^c$ tends to infinity to cover the whole parameter space. While the stepsize $h$ is just a mathematical artifact, the choice of $h$ actually influences the actual algorithm in Section \ref {sMMC} through the choice of the acceptance rate and the proposal probability in the following sections.

} } 

\subsection {The Boltzmann regime} \label {ssboltz}
In this section we derive the limiting equation for a large number $N$ of iterations in the algorithm (\ref {genmmc})
in the regime, where the proposal distribution $\tau (x_p|x,\kappa )$ is arbitrary, but the acceptance rate $\alpha $ is relatively small.
We rescale the acceptance rate to be of order $h=N^{c-1},\ 0<c<1$.
With this rescaling we have the following result:
{\proposition {If the acceptance rate $\alpha $ is uniformly of order $O(h)=O(N^{c-1})$, then  the solution $f(x,\kappa  ,s)$ will, for $h\rightarrow 0,\ N\rightarrow \infty $ converge against the solution of the kinetic integro - differential equation equation 
\begin {align} \label {lblzeq}
& \partial _s f(x,\kappa ,s)= \\
& -div _\kappa [\omega _1(x,\kappa ,s) f]
+
\int 
K(x,y,\kappa )f(y,\kappa ,s) \ dy
-
\int 
K(y,x,\kappa )f(x,\kappa ,s) \ dy \ ,
\end{align}
with the integral kernel $K$ given by 
\begin {equation}\label {lblzker}
K(x,y,\kappa )=\alpha (x,\kappa ,y,\kappa )\tau (x|y,\kappa ) \ .
\end {equation}
  }\label {prpcntindx}
} 
The proof of Proposition \ref {prpcntindx} is deferred to the Appendix.
\par
We now compute the limiting solution $f(x,\kappa ,\infty )$ of the solution of the kinetic transport equation (\ref {lblzeq}).
To this end we use the concept of entropy, i.e. we define a convex functional of the solution $f$ which will decay monotonically until a steady state is reached.
For the proof of convergence we will restrict ourselves to the case where the proposal and the acceptance rate do not depend on the moments $\kappa $, i.e. we set $\alpha =\alpha (x,y)$ and $\tau =\tau (x|y)$ in the equation (\ref {lblzeq}).
This allows to integrate the moment variable $\kappa $ out of equation (\ref {lblzeq}). So, under these assumptions, we consider the equation (\ref {lblzeq}) in its weak form (\ref {lblzwk})
\begin {equation} \label {lblzwk1}
\int \phi (x )\partial _s f(x ,s)\ dx =\int  [\phi (x)-\phi (y) ]
K(x,y)f(y,s) \ dxy \ ,
\end {equation}
with $K(x,y)=\alpha (x,y)\tau (x|y)$ and where we have chosen test functions $\phi $ which only depend on the state $x$.
 
We  define the limiting solution as a symmetrizer of the kernel $K$ in (\ref {lblzker}).
This is in kinetic theory sometimes called the concept of detail balance 
\cite {db1,db2}.
We formulate the:
\par 
{\bf Detail balance condition:} Let the the function $f_\infty (x)$ be defined by the symmetrizing condition
\begin {equation} \label {ldbblz}
\alpha (x,y)\tau (x|y)f_\infty (y)=\alpha (y,x)\tau (y|x)f_\infty (x)\ \forall \ x,y 
\end {equation} 

{\proposition {In the case that the acceptance probability $\alpha $, the proposal distribution $\tau $, and therefore the state probability $f(x,s)$, do no depend on the moments $\kappa $,
the probability distribution $f(x,s)$ will converge to the limiting distribution $f_\infty (x)$ for $s, N \rightarrow \infty $, which is given by the symmetry condition (\ref {ldbblz}).
  }\label {prpcnv} 
}
\par
The proof of Proposition \ref {prpcnv} is deferred to the Appendix.

\subsection {The Brownian motion regime  }\label {ssbrown}
In this section we derive an alternative continuum limit, based on a different limiting regime,  for a kinetic equation based on the evolution equation \eqref{lver4}.
This regime is based on the idea of Brownian motion and yields not an integral equation as in Section (\ref {ssboltz}), but a Fokker - Planck - like differential equation with a diffusive term.
According to \eqref{lver4} the discrete evolution  of the density $f$ is given in weak form by
\begin{align} \label {ldscbm}
& \int \phi (x,\kappa )f_{n+1}(x,\kappa )\ dx\kappa =
\phi (x,\omega _n(x,r))\alpha (x,\omega _n(x,r),y,\omega _n(y,r))\tau (x|y,r)f_n(y,r)\ dyrx
\\ 
& +\int \phi (y,\omega _n(y,r) )[1-\alpha (x,\omega _n(x,r),y,\omega _n(y,r))]\tau (x|y,r)f_n(y,r)\ dxry
\end{align}  
with the update of the moments given by
$$
\omega _n(x,\kappa )=\frac {\Delta \kappa (x)+n\kappa }{1+n}\ .
$$
Again, we replace this by a model with a continuous time step $s=nh,\ h=\frac  1 N$ with $N$ the total number of steps taken.
So $N\rightarrow \infty ,\ h\rightarrow 0$ holds in the continuum limit.
In Section (\ref {ssboltz}) we assume a small acceptance rate $\alpha $, proportional to the time step,
 and a proposed state $x$ with an arbitrary distance to the current state $y$.
Conversely,
the idea of the Brownian motion regime  is that the state $x$ is only increased incrementally,
that is with a mean proportional to the time step, but, other than in the Boltzmann regime in Section {\ref {ssboltz},
 the 
acceptance probability $\alpha $ can be chosen arbitrarily between zero and one.
In classical Brownian motion theory the proposed increment to the state is always accepted,yielding essentially a fractal for the particle path \cite {bm2,bm1}.
We modify the Brownian motion approach by including a variable acceptance rate $\alpha $, depending on the proposed  and the current state as well as the resulting higher order moments of the distribution.
Again we replace (\ref {ldscbm}) by a model with a continuous index $s=nh$
\begin{align} \label {lcntbm}
& \int \phi (x,\kappa )f(x,\kappa ,s+h)\ dx\kappa =
\phi (x,\omega (x,r,s))\alpha (x,\omega (x,r,s),y,\omega (y,r,s))\tau (x|y,r)f(y,r,s)\ dyrx\\
&+\int \phi (y,\omega (y,r,s) )[1-\alpha (x,\omega (x,r,s),y,\omega (y,r,s))]\tau (x|y,r)f(y,r,s)\ dxry,\
s=nh
\end{align}  
for all test functions $\phi(x,\kappa)    \in C^\infty_0( \mathbb{R}^{N_x} \times \mathbb{R}^{N_\kappa}; \mathbb{R}  )$, and with the update of the moments given by
\begin{equation} \label {lBMom}
\omega (x,\kappa ,s)=\frac  {h\Delta \kappa (x)+s\kappa }{h+s}\ .
\end{equation}  
For $\alpha =1$ (always accepting the proposal) this should yield the classical Fokker - Planck  equation 
obtained for the Brownian motion model. 
To separate the standard model from the influence of the acceptance rate,
we write this in terms of  a rejection rate $\beta =1-\alpha $. So, unconditional acceptance would mean $\beta =0$. This gives
\begin{align} \label {lbeta}
& \int \phi (x,\kappa )f(x,\kappa ,s+h)\ dx\kappa =
\int \phi (x,\omega (x,r,s))[1-\beta (x,\omega (x,r,s),y,\omega (y,r,s))]\tau (x|y,r)f(y,r,s)\ dyrx \\
& +\int \phi (y,\omega (y,r,s) )\beta (x,\omega (x,r,s),y,\omega (y,r,s))\tau (x|y,r)f(y,r,s)\ dxry,\
s=nh
\end{align}
Next, following the idea of the Brownian motion regime,
we write the proposal distribution in terms of increments to the current state,
setting $\tau (x|y,r)=\tau _1(x-y|y,r)$. So, $\tau _1$ is the distribution of the increment to the current state $y$.
Furthermore, we normalize the distribution $\tau _1$. We set 
$\tau _1(z|y,r) = \frac  {1}{(\sqrt h \sigma  (y,r))^d}\psi  (\frac  {z-hE(y,r)}{\sqrt h \sigma  (y,r)})$ with $\psi $ a normalized distribution with mean zero and the identity as covariance: $\int (1,z,zz^T)\psi (z)\ dz=(1,0,I)$. 
So, the expectation $hE(y,r)$ of $\tau _1$ is of the order of the stepsize  $h$ and  the standard deviation $\sqrt h \sigma  (y,r)$ is chosen such that it yields a diffusion term in the Brownian motion regime.
(Here $d$ denotes the dimension of the state vector $x$.)
So, altogether, we have
$$
\tau (x|y,r)=\frac  {1}{(\sqrt h \sigma  (y,r))^d}\psi (\frac  {x-y-hE(y,r)}{\sqrt h \sigma (y,r)})\ .
$$
With these assumptions (\ref {lbeta}) becomes
\begin{align} \label {lcntnrm}
& \int \phi (x,\kappa )f(x,\kappa ,s+h)\ dx\kappa = \\
& \int \phi (x,\omega (x,r,s))[1-\beta (x,\omega (x,r,s),y,\omega (y,r,s))] 
\frac  {1}{(\sqrt h \sigma (y,r))^d}\psi (\frac  {x-y-hE(y,r)}{\sqrt h \sigma (y,r)})f(y,r,s)\ dyrx
\\ 
&
+\int \phi (y,\omega (y,r,s) )\beta (x,\omega (x,r,s),y,\omega (y,r,s))
\frac  {1}{(\sqrt h \sigma (y,r))^d}\psi (\frac  {x-y-hE(y,r)}{\sqrt h \sigma (y,r)})f(y,r,s)\ dxry \ .
\end{align}

{\proposition {In the regime outlined above, i.e. the proposal distribution $\tau (x-y|y,\kappa )$ for the state increment 
has a mean $hE(y,\kappa )$ and a standard deviation $\sqrt h \sigma (y,\kappa )$,
the solution $f(x,\kappa ,s)$ will satisfy in the limit for many steps $N\rightarrow \infty ,h\rightarrow 0$ the Fokker - Planck equation

 \begin{align} \label {ldsf1}
& \partial _sf(x,\kappa ,s)=
-div _\kappa [\omega _1(x,\kappa ,s)f(x,\kappa ,s)]
\\
& +div_x[ E(x,\kappa ) (\beta (x,\kappa ,x,\kappa )-1)f(x,\kappa ,s)+\sigma  ^2\nabla _1\beta (x,\kappa ,x,\kappa ) f(x,\kappa ,s)]
\\
& +\frac 1 2 div _x \nabla _x [ \sigma  (x,\kappa )^2(1+\beta (x,\kappa ,x,\kappa ))f(x,\kappa ,s)]
\end{align} 
with the acceptance rate $\alpha =1-\beta $ ($\beta $ the rejection rate) and the increment $\omega _1(x,\kappa ,s)$ in the moments given by
$\omega _1(x,\kappa ,s)=\frac  {\Delta \kappa (x)-\kappa }{s}$.
}\label {prpbm}}
\par 
The proof of Proposition \ref {prpbm} is lengthy and deferred to the Appendix in Section {\ref {sappnd}.
\par 
{\remark {For a regime where we always accept the small increment in the state $x$ (for $\alpha (x,\kappa ,y,r)=1,\ \beta =0$)
this reduces to the classical Fokker - Planck equation for Brownian motion 
\cite {bm2}
with a drift term in the state and a drift term $\omega _1$ in the moments and a diffusive term caused by the variance  $\sigma ^2$ in the proposal increment. 
So, (\ref {ldsf1}) constitutes a modification to the classical Brownian motion model for variable rejection and  acceptance rates.
} }

\section {A micro--macro decomposition}\label {sec:micromacro}

The PDE \eqref{ldsf1} is possibly high--dimensional due to the dependence on $x$ and the vector of moments $\kappa.$ It is therefore prohibitively expensive to solve. In order to utilize the possibility to have probability density $\tau$ of the proposal depending on $x_n$ and moments $\kappa_n,$ we propose a micro--macro decomposition of the kinetic density $f.$ The microscopic part would then be solved using particles and an acceptance probability  $\alpha$ and proposal probability density $\tau$ depending solely on $x.$ The macroscopic part on the other hand would be solved with $\alpha,\tau$ solely depending on (a small number) of moments $\kappa.$ 
\par  
We consider the dynamics in the Boltzmann regime of Proposition \eqref{prpcntindx}. The case of the regime of the Brownian motion is similar. 
The probability density $f(x,\kappa,s)$ is the solution to the kinetic integro-diferential equation \eqref{lblzeq}:
\begin{align}
	\partial_s f(x,\kappa,s)       &=  \mathcal{Q}[f](x,\kappa,s), \\
	\mathcal{Q}[f](x,\kappa,s) &:= - div_\kappa  \omega_1(x,\kappa,s) f + \int K(x,y,\kappa) f(y,\kappa,s) - K(y,x,\kappa) f(x,\kappa,s) dy. 
\end{align}
The micro--macro decomposition seeks to find a solution of the form 
\begin{align}\label{decomp} 
	f(x,\kappa,s) = \zeta(s) f_{micro}(x,\kappa,s) + (1-\zeta(s)) f_{macro}(x,\kappa,s),
\end{align}
where $f_{micro}$ will denote the tail  and $f_{macro}$ will denote the bulk of the distribution of $f.$ Here, $ \zeta $ denotes the time--dependent splitting of the mass between tail and bulk. The update of this  quantity will be given below. The operator $\mathcal{Q}$ is linear in $f$ and therefore the decomposition \eqref{decomp} fulfills  for any $0< \gamma, \zeta < 1$ 
as

	\begin{align}
			\partial_s f_{micro}(x,\kappa,s) &= \mathcal{Q}[ \gamma f_{micro}(x,\kappa,s) + (1-\gamma) f_{macro}](x,\kappa,s), \\
			\partial_s f_{macro}(x,\kappa,s) &= \mathcal{Q}[ \frac{ \zeta (1-\gamma)}{1 - \zeta}   f_{micro}  (x,\kappa,s) + \frac{1 -2 \zeta + \zeta \gamma}{1-\zeta} f_{macro}](x,\kappa,s)  \\ 
			& - \frac{\zeta'(s)}{1-\zeta(s)} ( f_{micro} - f_{macro} )(x,\kappa,s). \label{eq-g}
		\end{align}
\\
The quantity $\zeta(s)$ distributes the mass between the microscopic $f_{micro}$ and the macroscopic $f_{macro}$ part of the distribution. The value of $\gamma$ is a parameter for the design of the method. Possible choices are for example $\gamma=1$ which leads to the system 
\begin{align}
	\partial_s f_{micro}(x,\kappa,s) &= \mathcal{Q}[  f_{micro}(x,\kappa,s) ], \\
	\partial_s f_{macro}(x,\kappa,s) &= \mathcal{Q}[f_{macro}](x,\kappa,s) - \frac{\zeta'(s)}{1-\zeta(s)} ( f_{micro} - f_{macro} )(x,\kappa,s).
\end{align}
while the choice $\gamma=\zeta$ leads to 
\begin{align}
	\partial_s f_{micro}(x,\kappa,s) = \mathcal{Q}[ \zeta f_{micro} + (1-\zeta) f_{macro}](x,\kappa,s), \\
	\partial_s f_{macro}(x,\kappa,s) = \mathcal{Q}[  \zeta  \;  f_{micro}   +(1-\zeta)  f_{macro}](x,\kappa,s) - \frac{\zeta'(s)}{1-\zeta(s)} ( f_{micro} - f_{macro} )(x,\kappa,s).
\end{align}

In particular, we solve the microscopic part $f_{micro}$ using the Metropolis Monte--Carlo method without dependence on $\kappa,$ i.e., $f_{micro}=f_{micro}(x,s),$ while the macroscopic part of the distribution $f_{macro}$ is solved by updating only moments in  $\kappa.$  In this case, we may assume that $\alpha, \tau$ are independent on $x$ and there exists a solution 
\begin{align}
	f_{macro}(x,\kappa,s) = \mathcal{M}_{\kappa(s)}(x) 
\end{align}
to equation \eqref{lblzeq}. The evolution $\kappa=\kappa(s)$ is then obtained by integration of equation \eqref{eq-g} against $(x,x^2, \dots)$ leading to 
\begin{align*}
	\int (x,x^2,\dots ) \partial_s  \mathcal{M}_{\kappa(s)}(x)  dx=
	\int (x,x^2, \dots )  \mathcal{Q}[ \frac{ \zeta (1-\gamma)}{1 - \zeta}   f_{micro}  (x,\kappa,s) + \frac{1 -2 \zeta + \zeta \gamma}{1-\zeta} \mathcal{M}_{\kappa(s)}] (x,\kappa,s) dx \\ 
	- \int (x,x^2, \dots ) \frac{\zeta'}{1-\zeta} ( f_{micro} - \mathcal{M}_{\kappa(s)} )(x,\kappa,s). 
\end{align*}
The later yields an update formula  for the moments of $f_{macro}.$ It remains to define the dynamics for $\zeta.$ The value of $\zeta(s+\zeta s)$ will be chosen to balance the variance of the data $\sigma^2(\nu_s)$ up to time $s$ and the variance of microscopic distribution $f_{micro}(s+\zeta s)$ and the bulk $f_{macro}(s).$  The details are given in the time-discrete case by equation \eqref{rate-r} leading to the definition of the distribution $\zeta(s).$

\section {Application to an Inverse Problem for the Lorentz System}\label {slorentznum}

The  example is taken from \cite[Example 2.6]{MR3363508}: the Lorenz '63 is a continuous dynamical  system that is known to exhibit sensitivity to initial conditions $v_0 \in \mathbb{R}^3$ as well as to parameters $x=(a,b,c) \in \mathbb{R}^3$.  Let $v=(v_1,v_2,v_3)(\cdot) \in C^1(0,T;\mathbb{R}^3)$ be the solution to the  set of ordinary differential equation up to a given time $T>0$ for a given set of parameters $x.$ 
\begin{subequations}\label{v63} \begin{align}
	\frac{d}{dt} v_1(t) &= a(v_2(t)-v_1(t)), \; v_1(0)=v_{1,0},  \\
	\frac{d}{dt} v_2(t) &= - a  v_1(t) - v_2(t) - v_1(t) v_3(t), \; v_2(0)=v_{2,0}, \\
	\frac{d}{dt} v_3(t) &= v_1(t) v_2(t) - b v_3(t) - b ( c + a ), \; v_3(0)=v_{3,0}.	
\end{align}
As in \cite{MR3363508}, the initial data is given by 
\begin{align}
	v_{i,0}=1, i=1,\dots,3,
\end{align}
\end{subequations}
We denote by $v(t;x)$  the solution to  system \eqref{v63} at time $t \in [0,T]$ for  given parameters $x.$ 
\par 
The problem is now rewritten in the form of problem \eqref{lminx}. In the notation For some fixed time $t$, the model $\mu(x)$ is hence given by 
\begin{align*}
	\mu(x) = v(t;x): \mathbb{R}^3 \to \mathbb{R}^3. 
\end{align*}
We consider two different scenarios for the definition of the set $\nu$ of data points.  In the following $K_{\max}$ is the number of data points and $x^* =\left(10, \frac83, 28 \right).$ 

{\remark {The system (\ref {v63}) is chaotic.
This makes a more detailed estimate of the parameters via the iterative Bayesian approach necessary \cite{MR3363508}.
}}

\subsection{Fixed terminal time}
In the first example, the set of observed data  is obtained as solution to \eqref{v63} at time $T$ using randomly perturbed parameters, i.e., 
$	\nu = \{z_i: i=1,\dots, K_{\max} \},$ and 
\begin{align}
	z_i = v(T, x^* + \xi_i ), i=1,\dots,K_{\max}.
\end{align}
Here, $\xi_i$ are the realization of a uniformly distributed random variable 
\begin{align} \label{v63xi}\xi \sim \otimes_{j=1}^3 \mathcal{U}(- a^j, a^j  ) .\end{align} 
We consider the Metropolis  Monte Carlo method {\em without }  information of the moments $\kappa.$  In the simulations, we compare using  a Gaussian distribution or a gradient based approach   to provide the proposal point $x_p.$  The further details of the implementation are stated in Section~\ref{details}.  In Figure \ref{fig-example1}  and Figure \ref{fig-example2} we compare the two strategies for generating proposals. While the true mean  is not exactly recovered, we observe that both strategies succeed in solving the minimization problem \eqref{lminx}. The histograms show a concentration close to the true mean. The later is obtained as $v(T,x^*)$ where $x^*$ is the true set of parameters without perturbation.

\begin{figure}[htb]\label{fig-example1}
	\includegraphics[width=.45\textwidth]{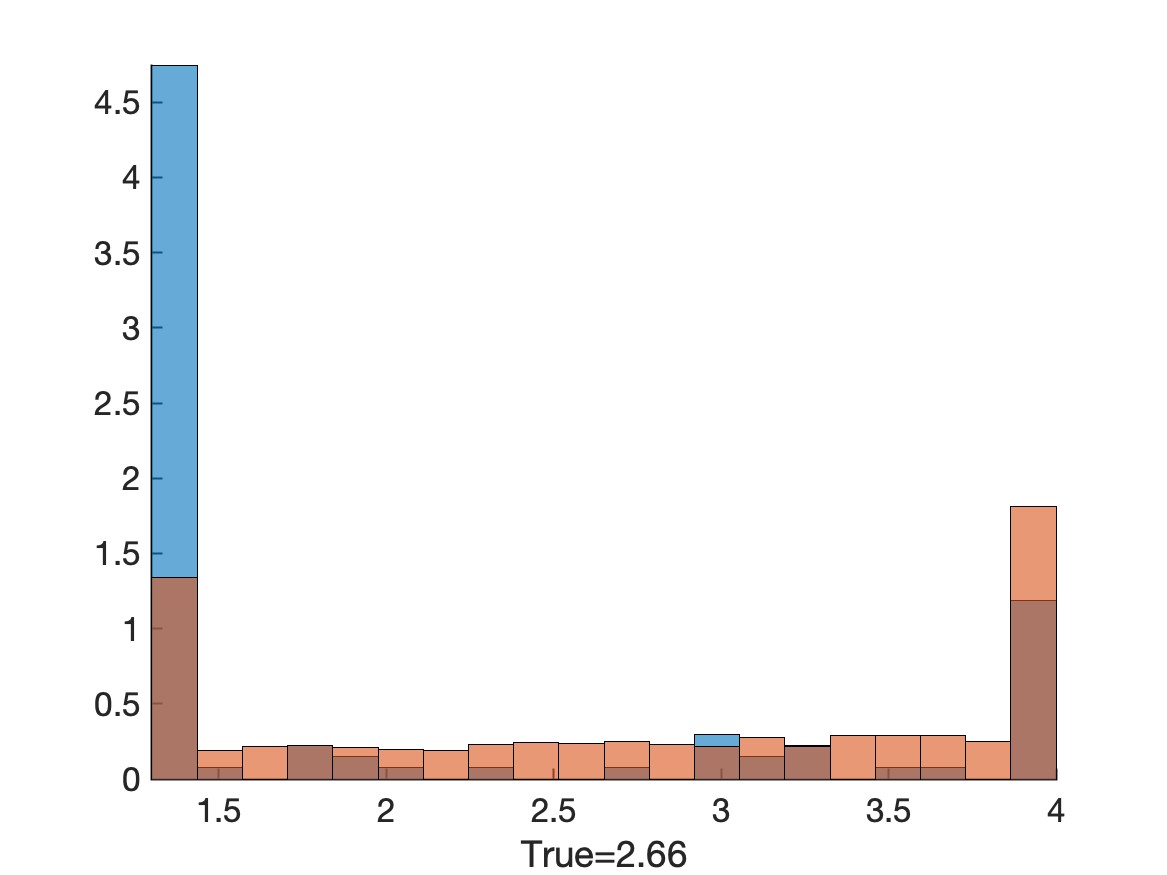}
\includegraphics[width=.45\textwidth]{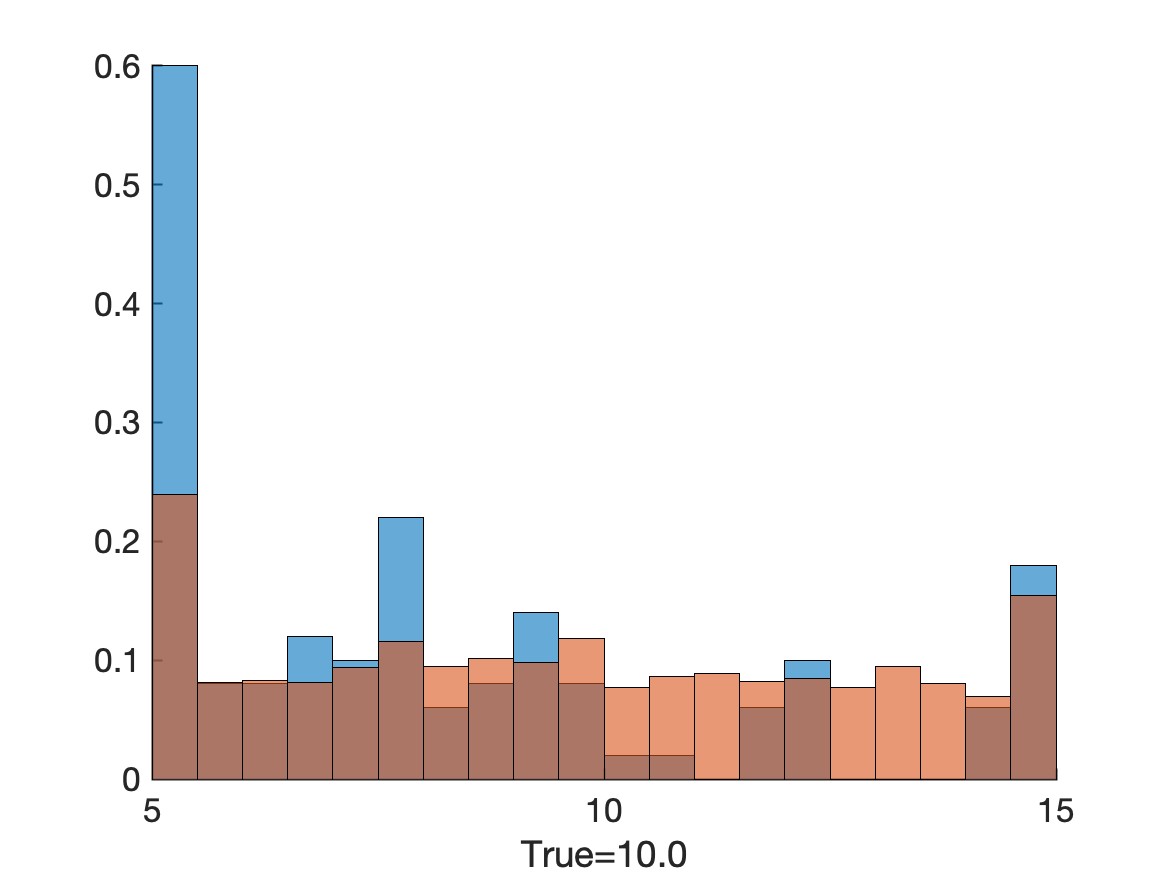}
	\\
	\includegraphics[width=.45\textwidth]{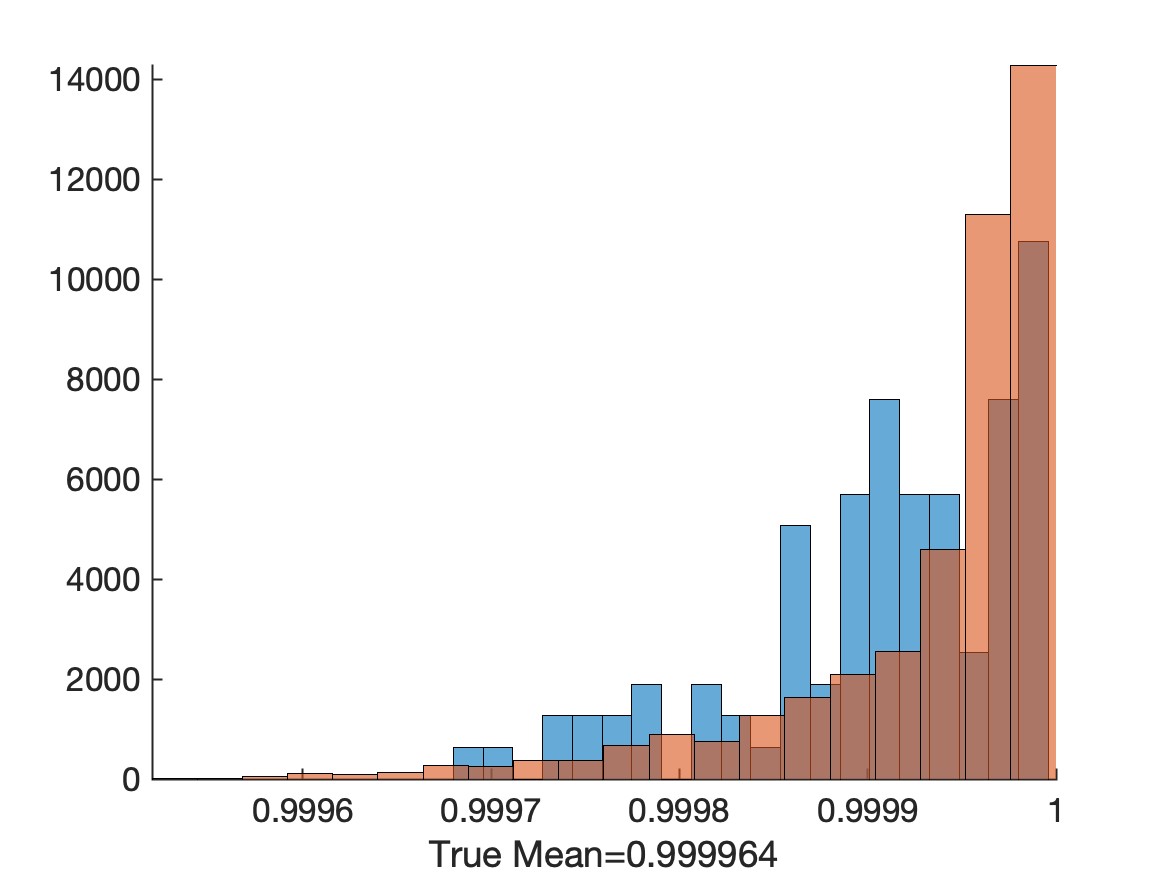}
	\includegraphics[width=.45\textwidth]{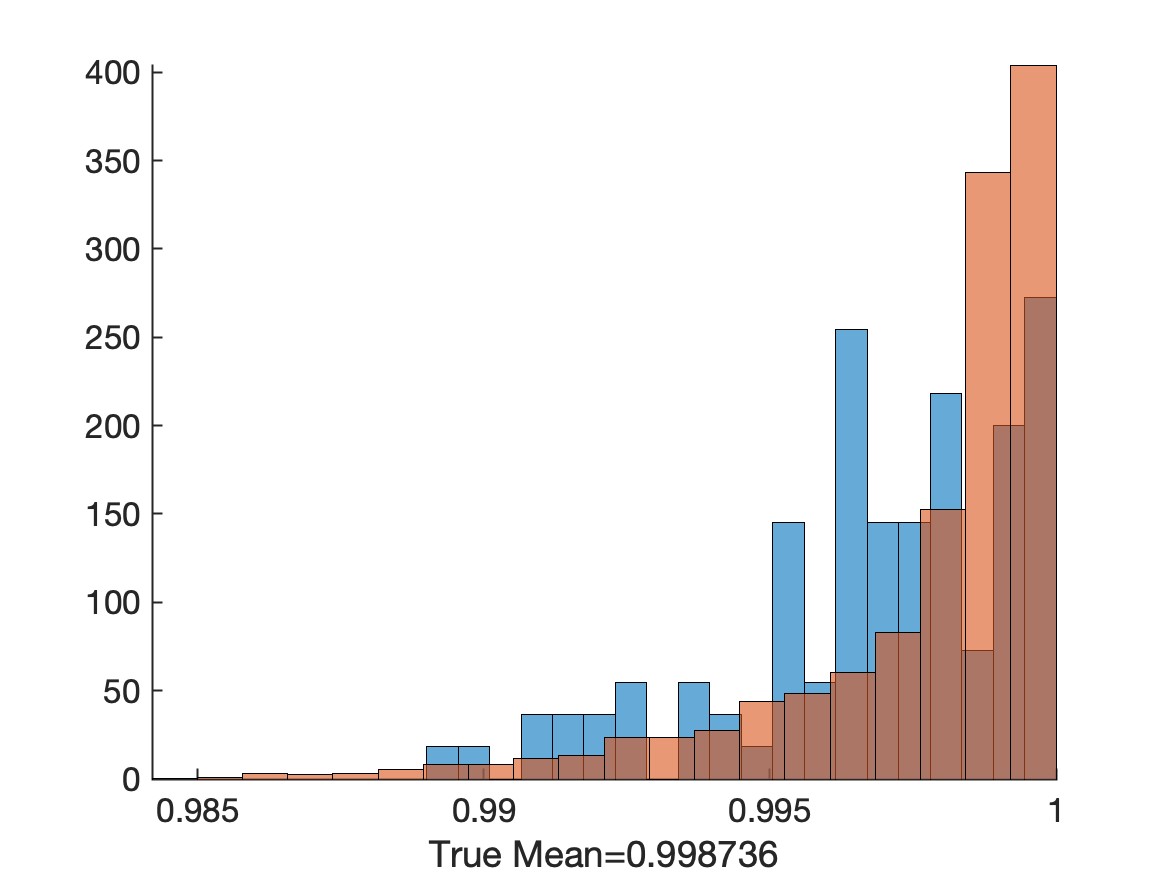}

	\caption{Fixed terminal time simulated by Metropolis Monte--Carlo using Gaussian proposals. 
		 The initial distribution $P_0$ is depicted in blue as a histogram. The histogram of the terminal distribution
		 $P_N$ is shown in red. In the top part we show the histogram of the parameters, in the bottom part the corresponding histogram of the model evaluations, i.e. $\mu$. With the parameters and model evaluations similar, we show $x_2$ and $x_3$ and $v_2$ and $v_3$ in the lower part of the diagram, respectively. The value 'true mean' represents the solution $v(T,x^*)$ for the optimal parameter $x^*.$ }
	\end{figure}

\begin{figure}[htb]\label{fig-example2}
	\includegraphics[width=.45\textwidth]{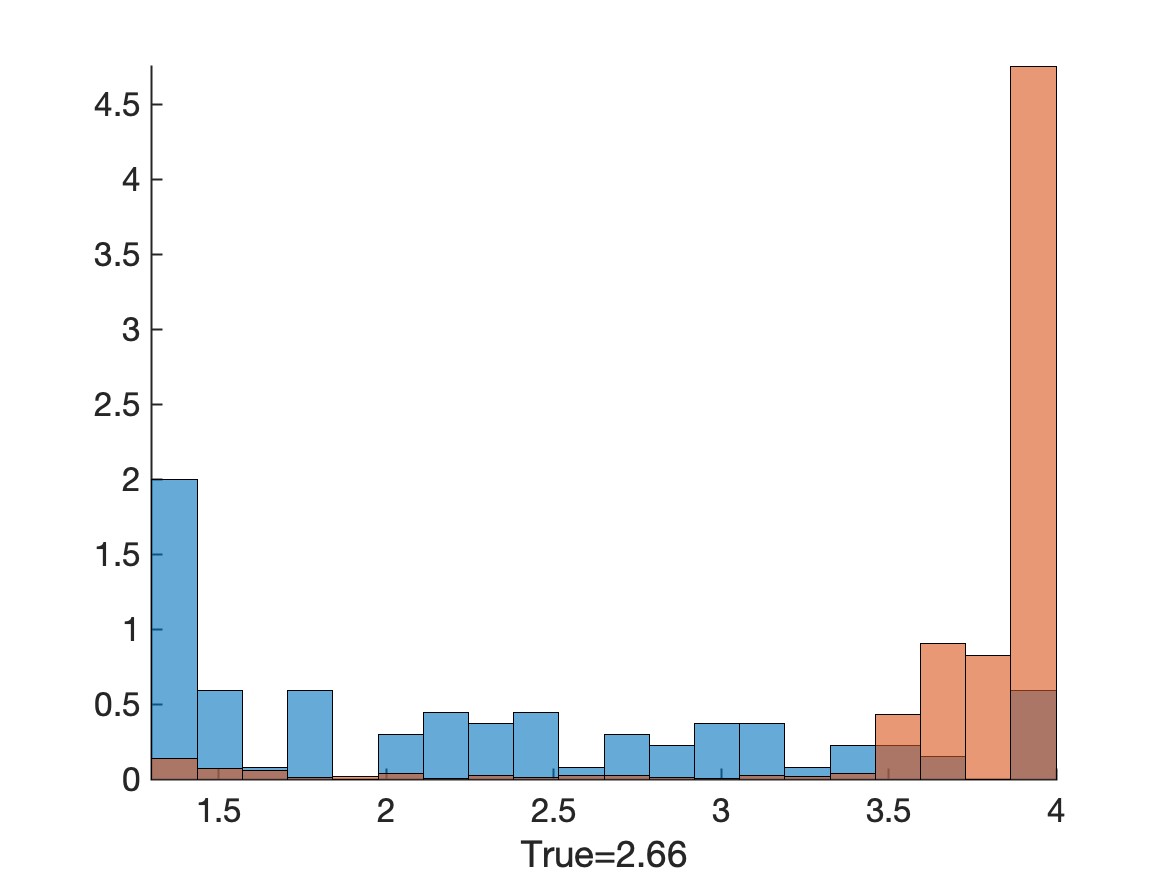}
	\includegraphics[width=.45\textwidth]{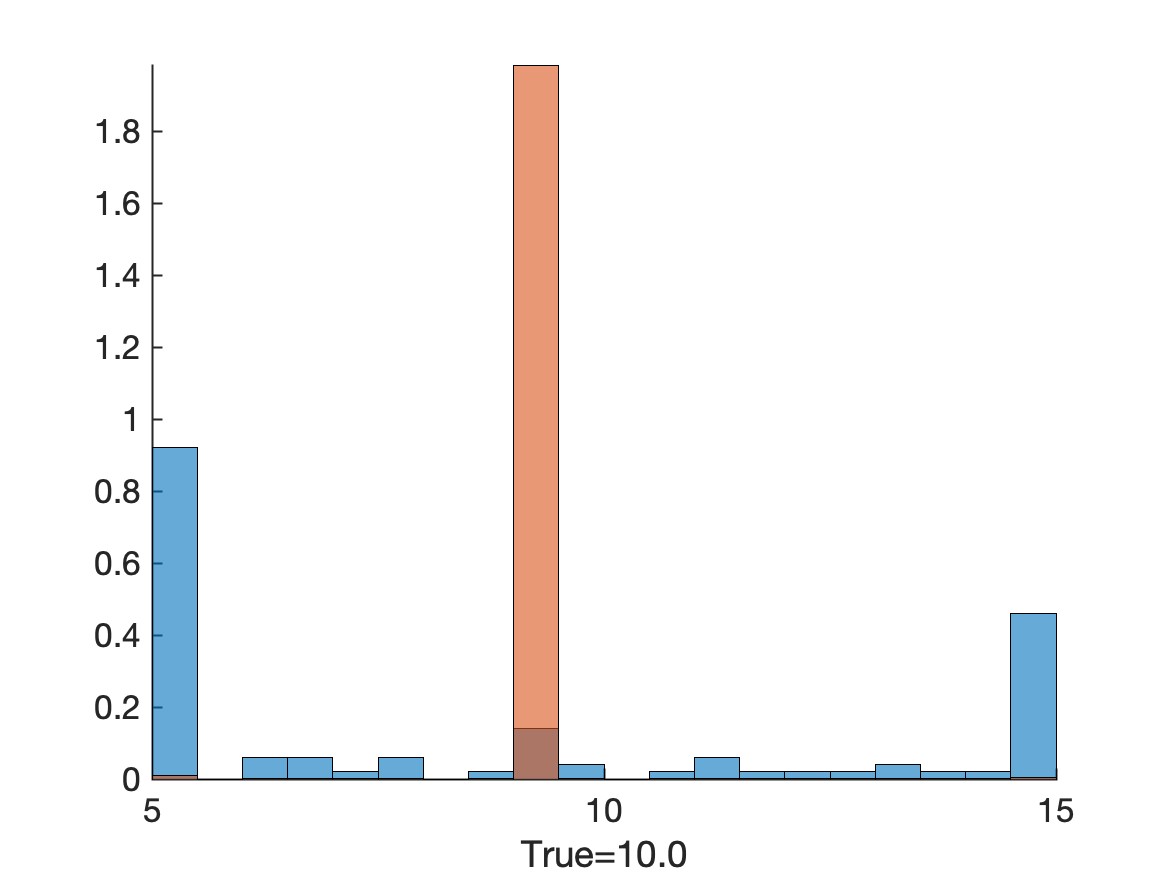}
	\\
	\includegraphics[width=.45\textwidth]{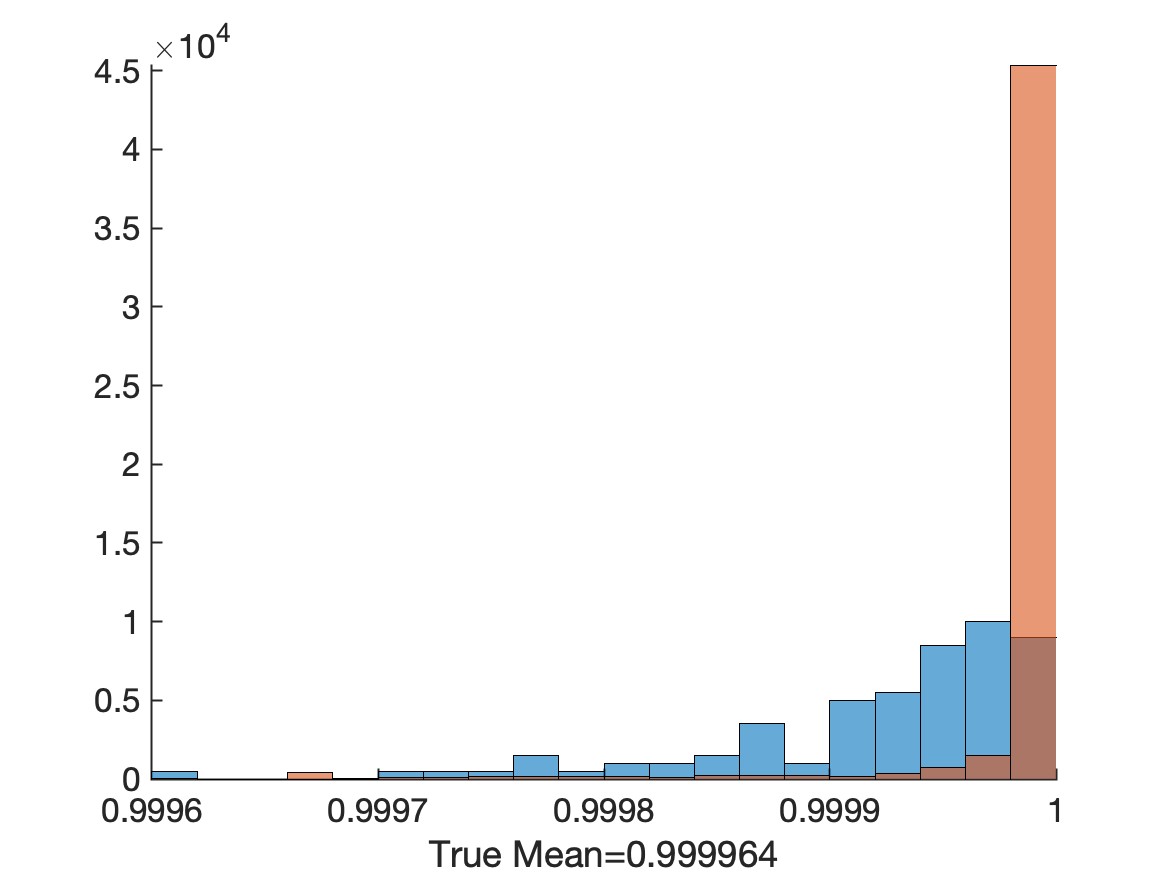}
	\includegraphics[width=.45\textwidth]{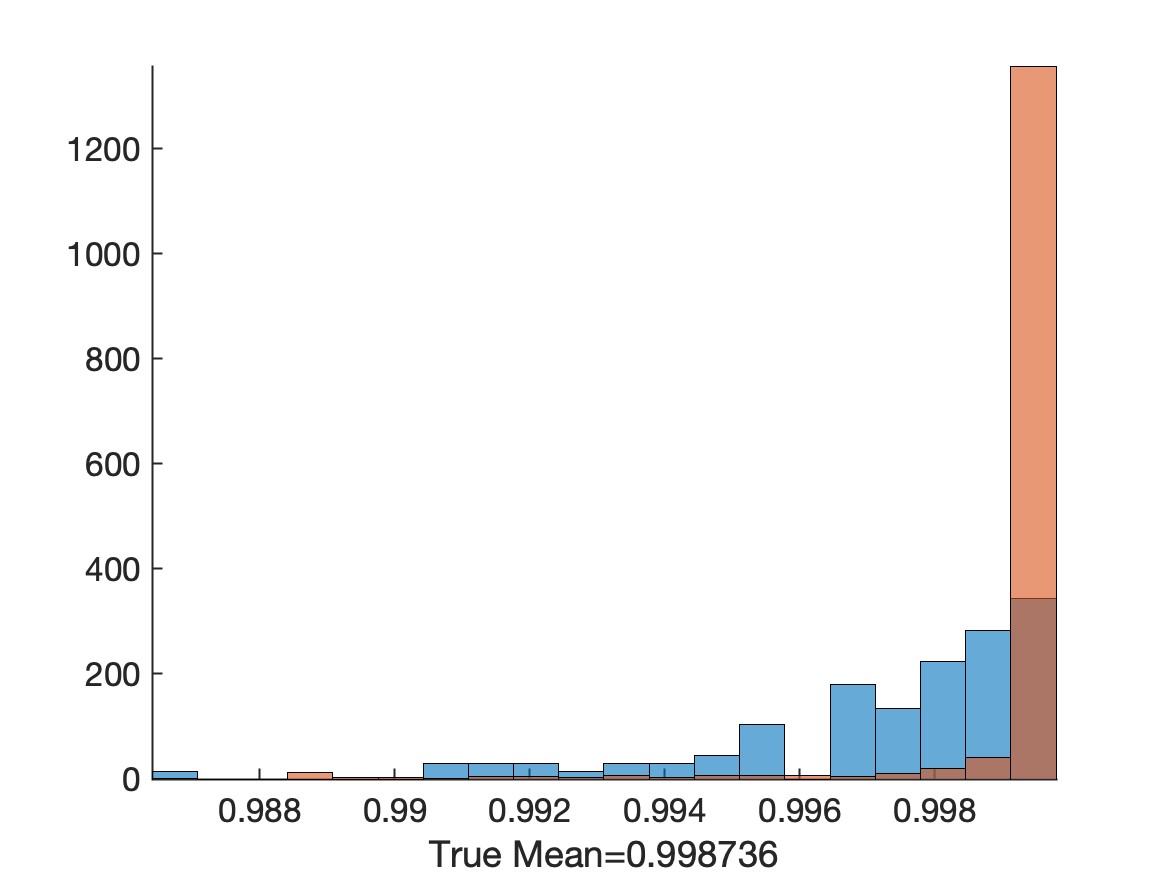}
	
	\caption{Fixed terminal time simulated using gradient--based proposal updates.  The initial distribution $P_0$ is depicted in blue as a histogram. The histogram of the terminal distribution
		$P_N$ is shown in red. In the top part we show the histogram of the parameters, in the bottom part the corresponding histogram of the model evaluations, i.e. $\mu$. With the parameters and model evaluations similar, we show $x_2$ and $x_3$ and $v_2$ and $v_3$ in the lower part of the diagram, respectively. The value 'true mean' represents the solution $v(T,x^*)$ for the optimal parameter $x^*.$ }
\end{figure}

\subsection{Running terminal time}
In the second example, we consider a problem where the time horizon is not fixed. Denote by $t_i =  i \Delta t$ for $i=1,\dots,K_{\max}$ and with $K_{\max} \Delta t = T$ for a fixed  terminal time $T.$   The set of data points $\nu = \{z_i: i=1,\dots, K_{\max} \}$ is  given by the solution at time $t_i$ for perturbed parameters: 
\begin{align}\label{dataset-ex}
	z_i = v(t_i, x^* + \xi_i )
\end{align}
where $\xi_i$ are realizations of a random variable $\xi$ distributed as in equation \eqref{v63xi}.  We compare the Metropolis Monte Carlo method without and with using moment information  on the identification problem \eqref{lminx}. The proposal $x_p$ is obtained using a micro-macro decomposition using an indicator $\zeta$ as outlined in Section \ref{sec:micromacro}. The further details of the implementation are stated in Section~\ref{details}. The results are depicted in Figure \ref{fig-example3} and Figure \ref{fig-example4}, respectively. We observe that the micro--macro decomposition is feasible and leads to similar histograms compared with the Metropolis Monte Carlo method. The rate of the macroscopic updates is about $1\%$ in the reported results.

\begin{figure}[htb]\label{fig-example3}
	\includegraphics[width=.45\textwidth]{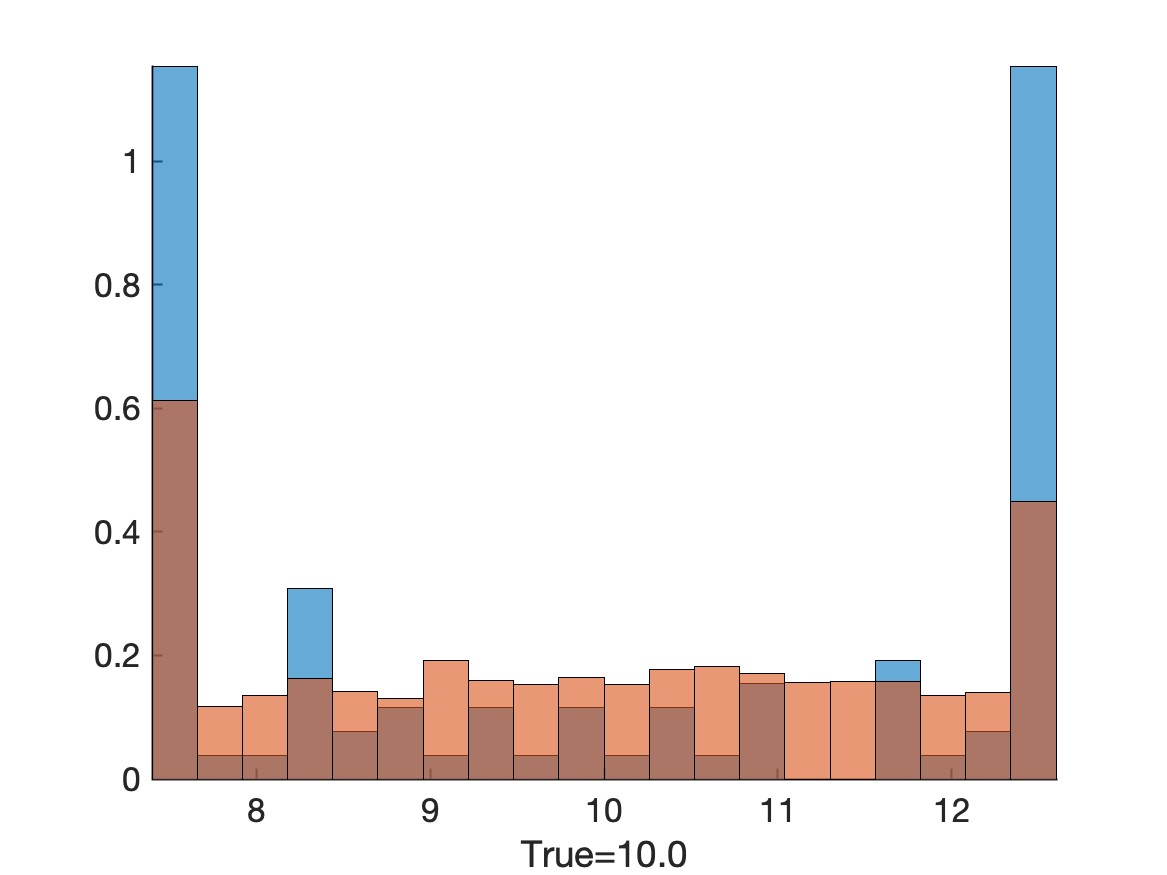}
	\includegraphics[width=.45\textwidth]{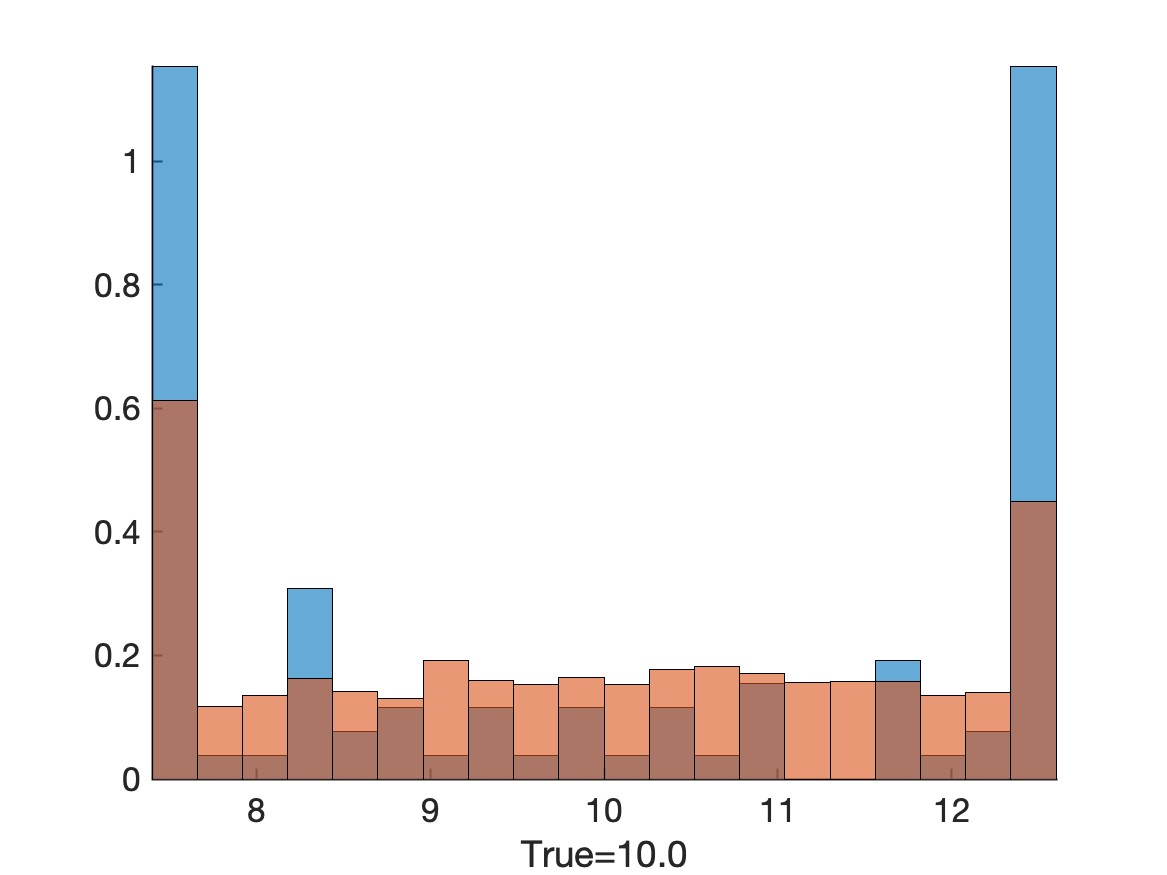}
	\\
	\includegraphics[width=.45\textwidth]{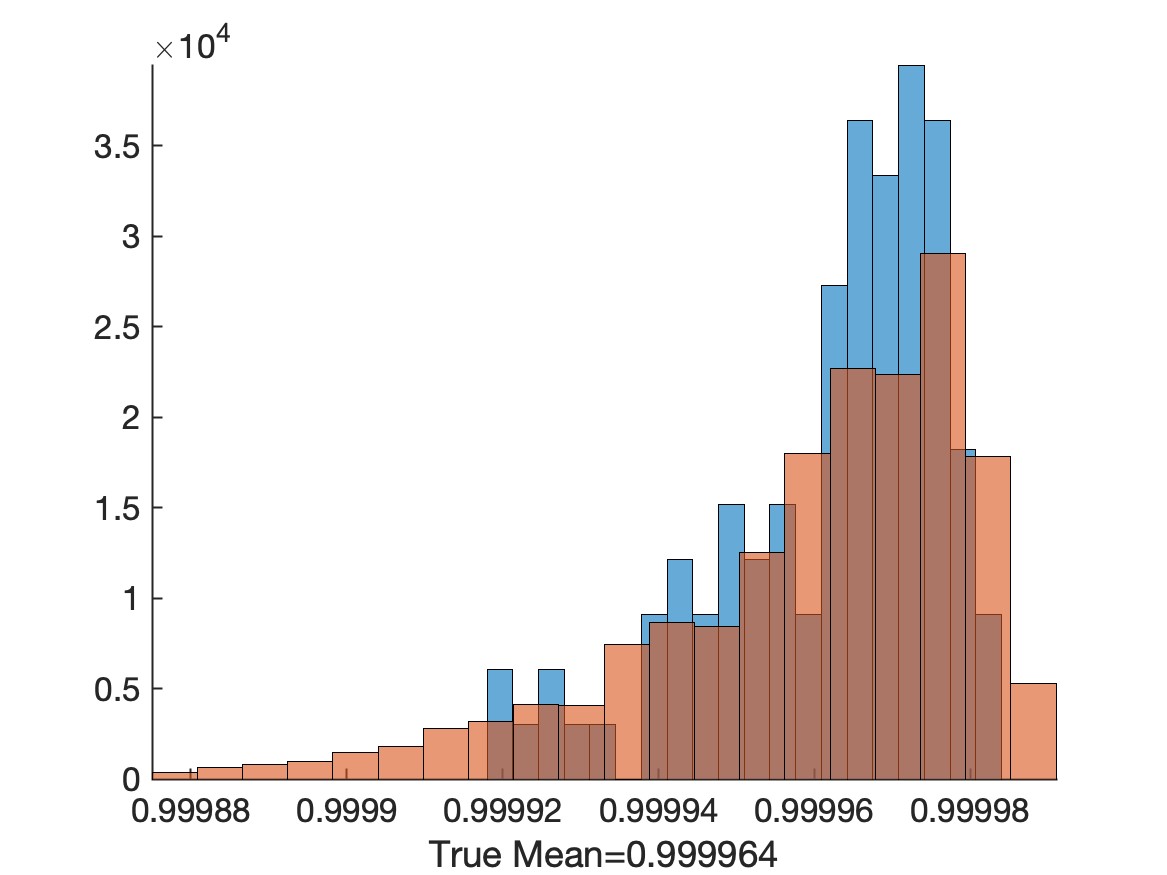}
	\includegraphics[width=.45\textwidth]{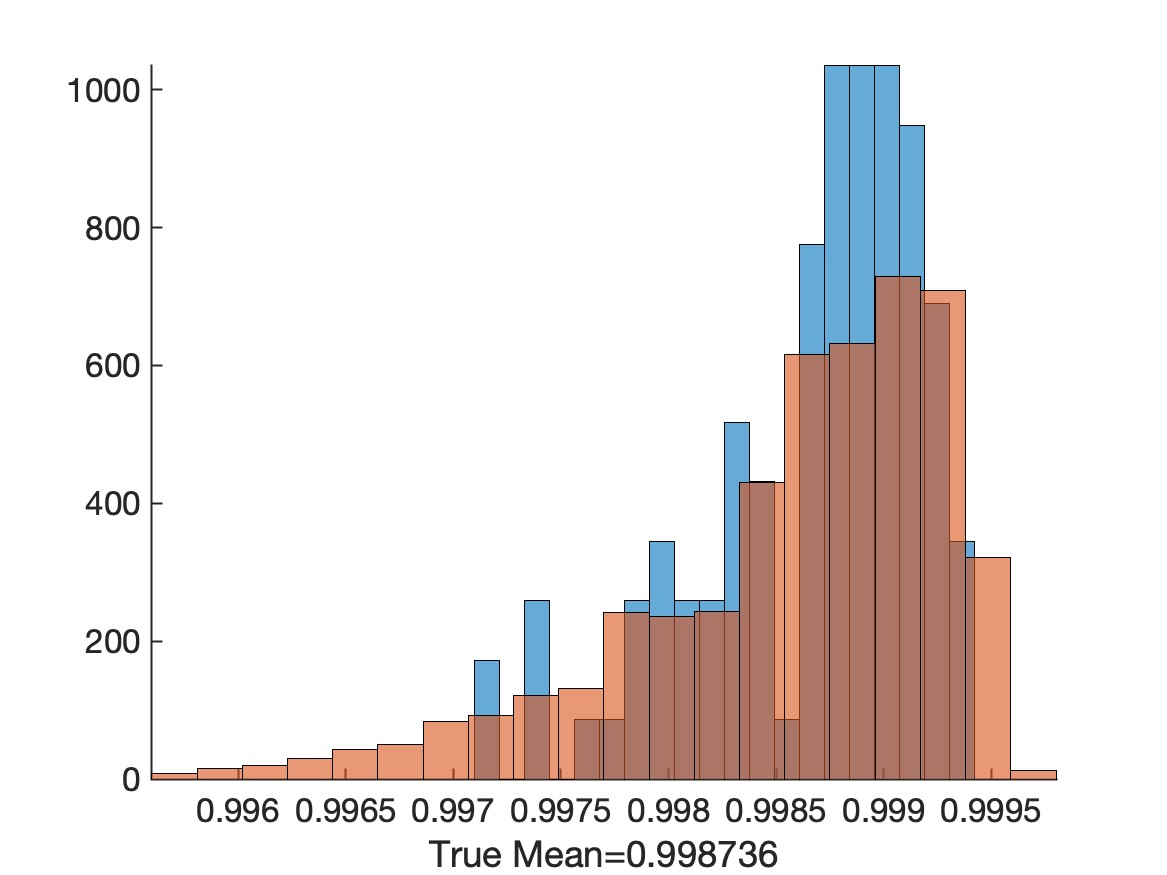}
	
	\caption{Running terminal time simulated using Metropolis Monte Carlo method.  The initial distribution $P_0$ is depicted in blue as a histogram. The histogram of the terminal distribution
		$P_N$ is shown in red. In the top part we show the histogram of the parameters, in the bottom part the corresponding histogram of the model evaluations, i.e. $\mu$. With the parameters and model evaluations similar, we show $x_2$ and $x_3$ and $v_2$ and $v_3$ in the lower part of the diagram, respectively. The value 'true mean' represents the solution $v(T,x^*)$ for the optimal parameter $x^*.$ }
\end{figure}

\begin{figure}[htb]\label{fig-example4}
	\includegraphics[width=.45\textwidth]{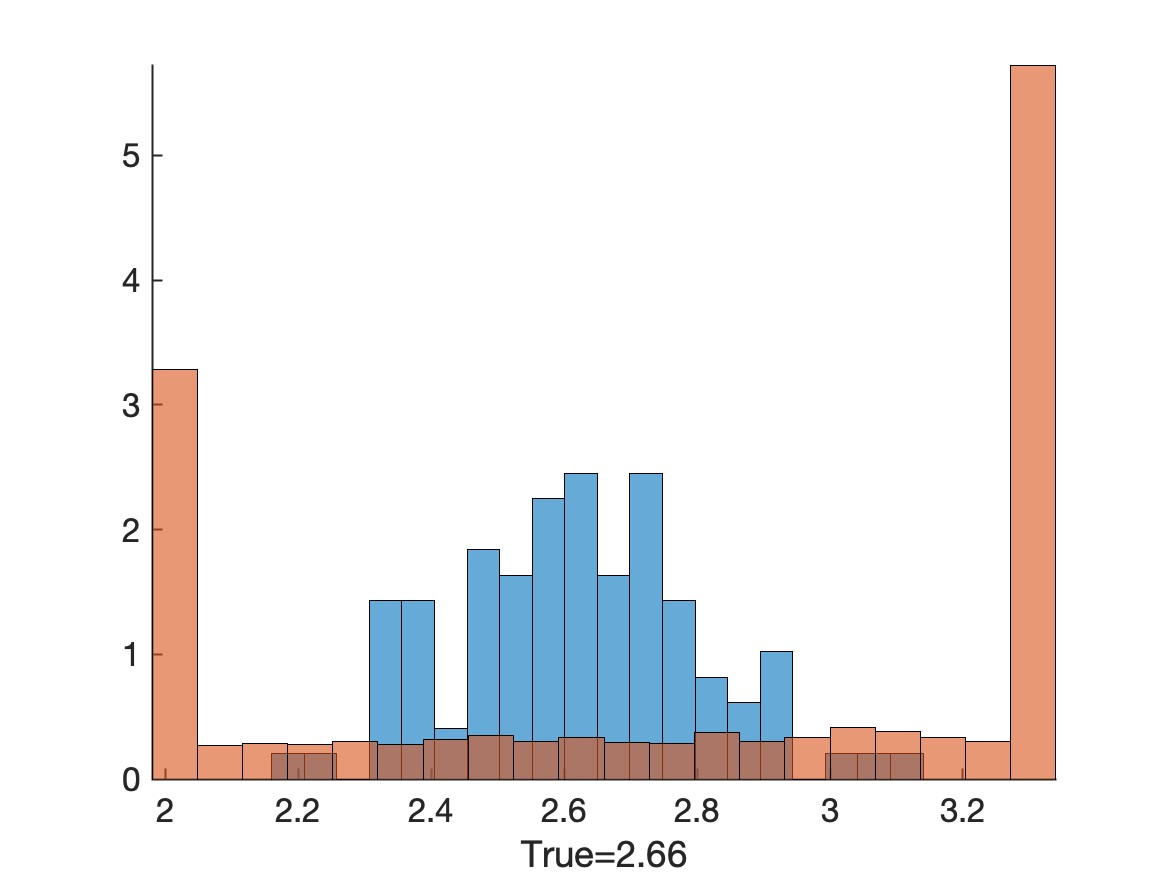}
	\includegraphics[width=.45\textwidth]{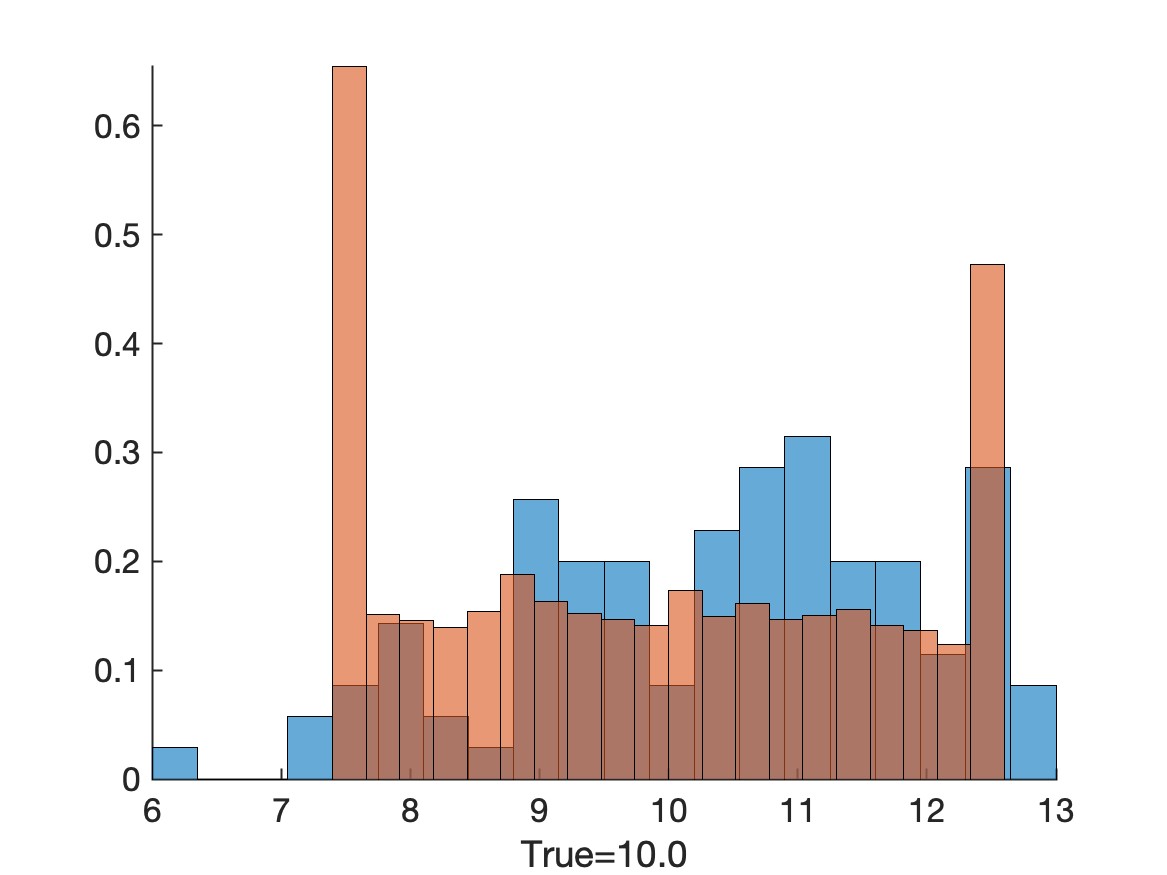}
	\\
	\includegraphics[width=.45\textwidth]{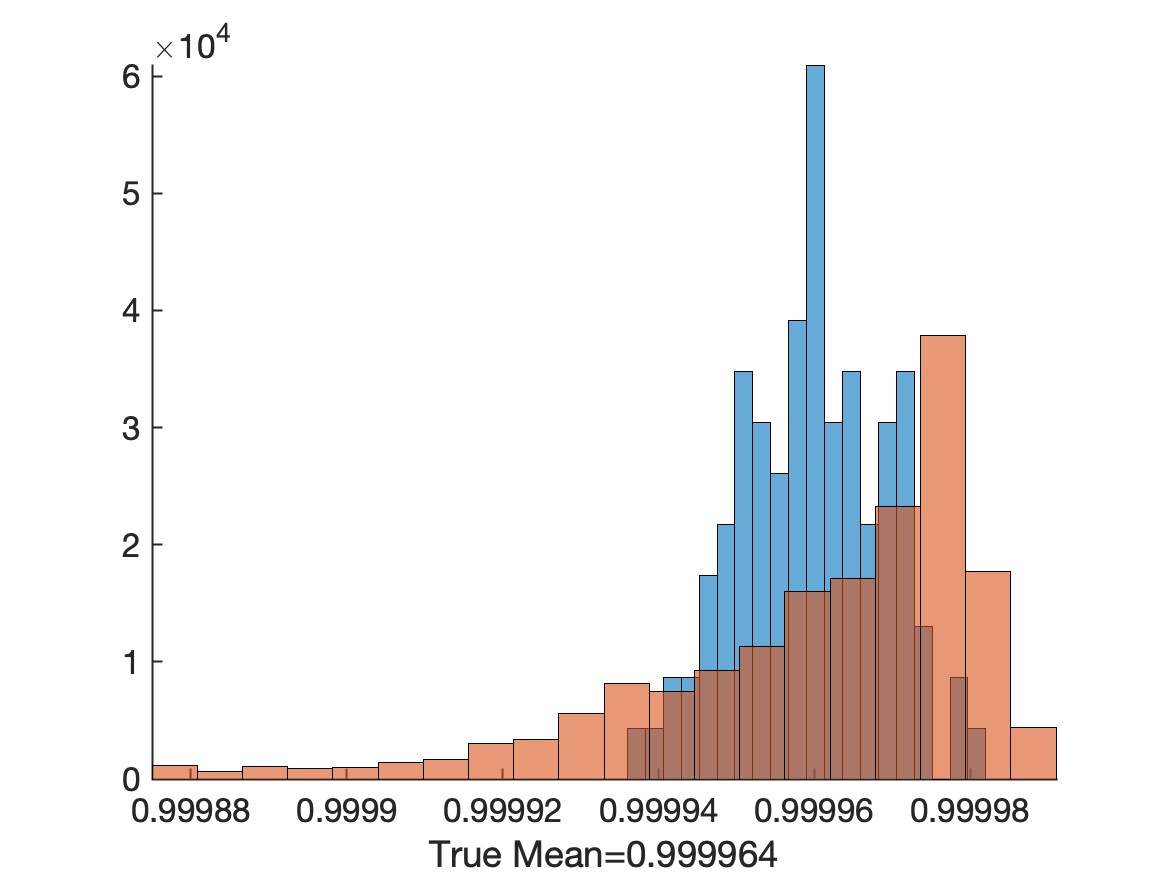}
	\includegraphics[width=.45\textwidth]{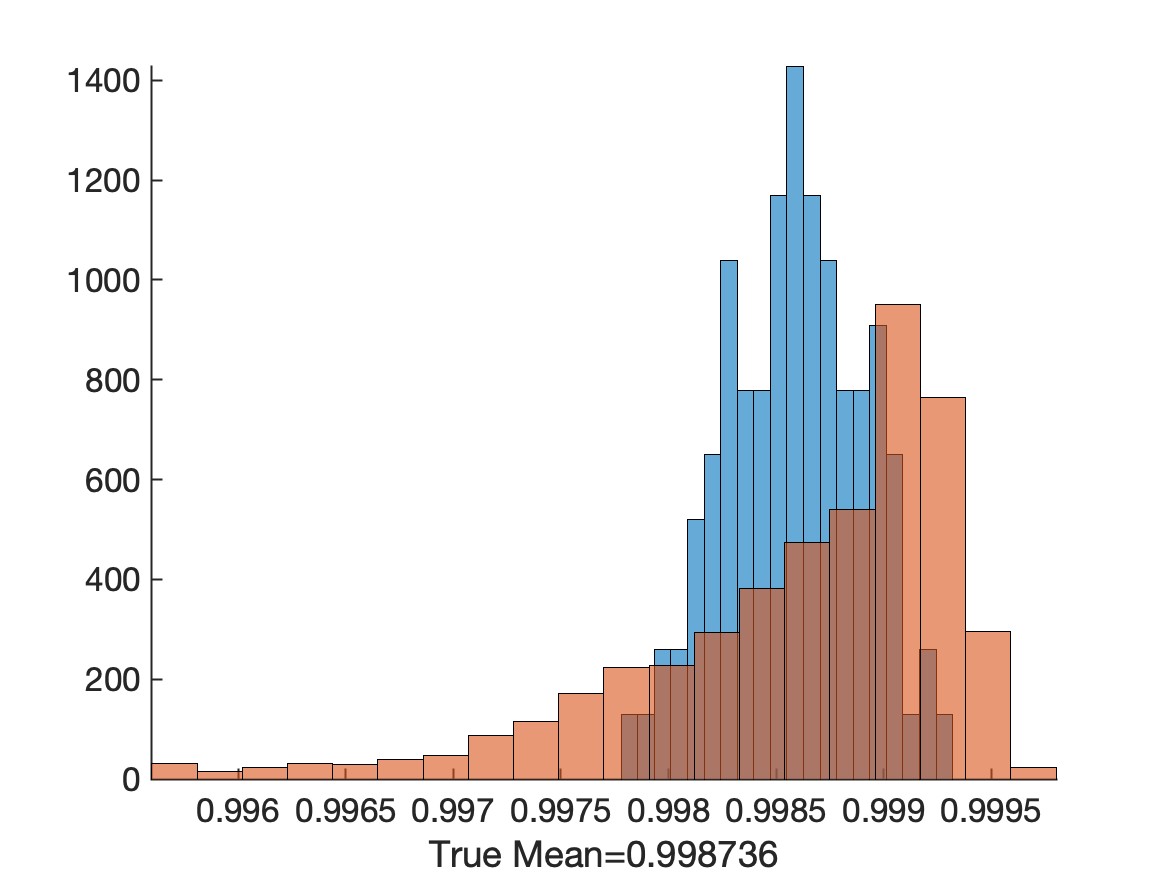}
	
	\caption{Running terminal time simulated using micro--macro decomposition.  The initial distribution $P_0$ is depicted in blue as a histogram. The histogram of the terminal distribution
		$P_N$ is shown in red. In the top part we show the histogram of the parameters, in the bottom part the corresponding histogram of the model evaluations, i.e. $\mu$. With the parameters and model evaluations similar, we show $x_2$ and $x_3$ and $v_2$ and $v_3$ in the lower part of the diagram, respectively. The value 'true mean' represents the solution $v(T,x^*)$ for the optimal parameter $x^*.$ }
\end{figure}


\section{Summary}
In this paper we developed a kinetic convergence theory for Metropolis Monte Carlo algorithms.
The application to Bayesian type inverse problems, where the result is not one optimal parameter but
a probability distribution in parameter space which optimally fits the given data, is considered.  
The kinetic theory allows for the theoretical reduction to lower dimensional problems which, in turn, allows for the efficient use of improved predictors an proposal distributions.
All this leads to  iterative Bayesian estimation procedures with a significant increase in the quality of the posterior distribution as well as an increase in computational efficiency.
\par 
\vspace{1cm}
\noindent {\small {\bf Acknowledgments}
The authors thank the Deutsche Forschungsgemeinschaft (DFG, German Research Foundation) for the financial support through 320021702/GRK2326,  333849990/IRTG-2379, B04, B05 and B06 of 442047500/SFB1481, HE5386/19-3,22-1,23-1,25-1,26-1,27-1  and support  from the European Unions Horizon Europe research and innovation programme under the Marie Sklodowska-Curie Doctoral Network Datahyking (Grant No. 101072546).

}
\bibliographystyle{plain} 
\bibliography{refs} 

\begin{thebibliography}{10}

\bibitem{db1}
M~Bauer and F~Cornu.
\newblock Local detailed balance: a microscopic derivation.
\newblock {\em Journal of Physics A: Mathematical and Theoretical},
  48(1):015008, dec 2014.

\bibitem{brooks2011handbook}
S.~Brooks, A.~Gelman, G.~Jones, and X.~Meng.
\newblock {\em Handbook of Markov Chain Monte Carlo}.
\newblock CRC press, 2011.

\bibitem{MR3126997}
Tan Bui-Thanh, Omar Ghattas, James Martin, and Georg Stadler.
\newblock A computational framework for infinite-dimensional {B}ayesian inverse
  problems {P}art {I}: {T}he linearized case, with application to global
  seismic inversion.
\newblock {\em SIAM J. Sci. Comput.}, 35(6):A2494--A2523, 2013.

\bibitem{Somersalo2007}
D.~Calvetti and E.~Somersalo.
\newblock {\em Introduction to {B}ayesian {S}cientific {C}omputing}.
\newblock Springer, 2007.

\bibitem{weather1}
Rafael Cano, Carmen Sordo, and Jos{\'e} Guti{\'e}rrez.
\newblock Applications of bayesian networks in meteorology.
\newblock {\em Wiley Interdisciplinary Reviews: Computational Statistics},
  2004.

\bibitem{db2}
Hudong Chen.
\newblock {$H$}-theorem and generalized semi-detailed balance condition for
  lattice gas systems.
\newblock {\em J. Statist. Phys.}, 81(1-2):347--359, 1995.

\bibitem{bm2}
Albert Einstein.
\newblock {\em Investigations on the theory of the {B}rownian movement}.
\newblock Dover Publications, Inc., New York, 1956.
\newblock Edited with notes by R. F\"{u}rth, Translated by A. D. Cowper.

\bibitem{MR1009037}
Heinz~W. Engl, Karl Kunisch, and Andreas Neubauer.
\newblock Convergence rates for {T}ikhonov regularisation of nonlinear
  ill-posed problems.
\newblock {\em Inverse Problems}, 5(4):523--540, 1989.

\bibitem{MR2260716}
Dani Gamerman and Hedibert~Freitas Lopes.
\newblock {\em Markov chain {M}onte {C}arlo}.
\newblock Texts in Statistical Science Series. Chapman \& Hall/CRC, Boca Raton,
  FL, second edition, 2006.
\newblock Stochastic simulation for Bayesian inference.

\bibitem{MR2858443}
Charles~J. Geyer.
\newblock Introduction to {M}arkov chain {M}onte {C}arlo.
\newblock In {\em Handbook of {M}arkov chain {M}onte {C}arlo}, Chapman \&
  Hall/CRC Handb. Mod. Stat. Methods, pages 3--48. CRC Press, Boca Raton, FL,
  2011.

\bibitem{MR2884617}
M.~Hairer, A.~Stuart, and J.~Voss.
\newblock Signal processing problems on function space: {B}ayesian formulation,
  stochastic {PDE}s and effective {MCMC} methods.
\newblock In {\em The {O}xford handbook of nonlinear filtering}, pages
  833--873. Oxford Univ. Press, Oxford, 2011.

\bibitem{MR4565586}
Bastian Harrach, Tim Jahn, and Roland Potthast.
\newblock Regularizing linear inverse problems under unknown non-{G}aussian
  white noise allowing repeated measurements.
\newblock {\em IMA J. Numer. Anal.}, 43(1):443--500, 2023.

\bibitem{MR3363437}
W.~K. Hastings.
\newblock Monte {C}arlo sampling methods using {M}arkov chains and their
  applications.
\newblock {\em Biometrika}, 57(1):97--109, 1970.

\bibitem{MR4687001}
P.~H\'{e}as, F.~C\'{e}rou, and M.~Rousset.
\newblock Chilled sampling for uncertainty quantification: a motivation from a
  meteorological inverse problem.
\newblock {\em Inverse Problems}, 40(2):Paper No. 025004, 38, 2024.

\bibitem{Kaipio2005}
J.~Kaipio and E.~Somersalo.
\newblock {\em {S}tatistical and {C}omputational {I}nverse {P}roblems}.
\newblock 0066-5452. Springer-Verlag New York, 2005.

\bibitem{MR3363508}
Kody Law, Andrew Stuart, and Konstantinos Zygalakis.
\newblock {\em Data assimilation}, volume~62 of {\em Texts in Applied
  Mathematics}.
\newblock Springer, Cham, 2015.
\newblock A mathematical introduction.

\bibitem{MR4196544}
Aaron Myers, Alexandre~H. Thi\'{e}ry, Kainan Wang, and Tan Bui-Thanh.
\newblock Sequential ensemble transform for {B}ayesian inverse problems.
\newblock {\em J. Comput. Phys.}, 427:Paper No. 110055, 21, 2021.

\bibitem{MR4604099}
Richard Nickl.
\newblock {\em Bayesian non-linear statistical inverse problems}.
\newblock Zurich Lectures in Advanced Mathematics. EMS Press, Berlin, [2023]
  \copyright 2023.

\bibitem{MR4619819}
Daniel Sanz-Alonso, Andrew Stuart, and Armeen Taeb.
\newblock {\em Inverse problems and data assimilation}, volume 107 of {\em
  London Mathematical Society Student Texts}.
\newblock Cambridge University Press, Cambridge, 2023.

\bibitem{bm1}
Lutz Schimansky-Geier and Thorsten P\"{o}schel, editors.
\newblock {\em Stochastic dynamics}, volume 484 of {\em Lecture Notes in
  Physics}.
\newblock Springer-Verlag, Berlin, 1997.

\bibitem{MR2652785}
A.~M. Stuart.
\newblock Inverse problems: a {B}ayesian perspective.
\newblock {\em Acta Numer.}, 19:451--559, 2010.

\bibitem{Sullivan2015}
T.~Sullivan.
\newblock {\em Introduction to {U}ncertainty {Q}uantification}.
\newblock Springer, 2015.

\bibitem{startup2}
Leila Taghizadeh, Ahmad Karimi, Benjamin Stadlbauer, Wolfgang~J. Weninger,
  Eugenijus Kaniusas, and Clemens Heitzinger.
\newblock Bayesian inversion for electrical-impedance tomography in medical
  imaging using the nonlinear poisson boltzmann equation.
\newblock {\em Computer Methods in Applied Mechanics and Engineering},
  365:112959, 2020.

\bibitem{startup1}
Leila Taghizadeh, Amirreza Khodadadian, and Clemens Heitzinger.
\newblock The optimal multilevel monte-carlo approximation of the stochastic
  drift diffusion poisson system.
\newblock {\em Computer Methods in Applied Mechanics and Engineering},
  318:739--761, 2017.

\end{thebibliography}

\section {Appendix}\label {sappnd}
\subsection{Technical proofs}

{\bf Proof of Proposition \ref {prpcntindx}:}
\\
Separating the terms in (\ref {lver4})into those containing the acceptance rate and the term independent of $\alpha  $ gives, using
the fact that $\int \tau (x|y,r)\ dx=1$ holds, 
\begin {align}\label {lver5}
& \int  \phi (x,\kappa )f(x,\kappa ,s+h)\ dx\kappa =
\int  \phi (y,r+h\omega _1(y,r,s)) f(y,r,s) \ dyr +
\\ 
& h\int  [\phi (x,r+h\omega _1(x,r,s))-\phi (y,r+h\omega _1(y,r,s)) ]\times 
\\ & \alpha (x,r+h\omega _1(x,r,s),y,r+h\omega _1(y,r,s))\tau (x|y,r)f(y,r,s) \ dxyr 
\end{align}
The first integral on the right hand side of (\ref {lver5}) gives after Taylor expansion in the stepsize $h$
\begin {align}\label {lblzdrft}
& \int  \phi (y,r+h\omega _1(y,r,s)) f(y,r,s) \ dyr=
\\ & 
\int  \phi (y,r) f(y,r,s) \ dyr
+
\int h \nabla _2\phi (y,r)\cdot \omega _1(y,r,s) f(y,r,s) \ dyr
+o(h) \ .
\end{align}
The second term on the right hand side of (\ref {lver5}) gives, up to order $o(h)$
\begin {align} \label {lblzcll}
& h\int  [\phi (x,r+h\omega _1(x,r,s))-\phi (y,r+h\omega _1(y,r,s)) ]\times 
\\ & \alpha (x,r+h\omega _1(x,r,s),y,r+h\omega _1(y,r,s))\tau (x|y,r)f(y,r,s) \ dxyr =
\\ & h\int  [\phi (x,r)-\phi (y,r) ]
K(x,y,r)f(y,r,s) \ dxyr +o(h) \ ,
\end{align}
with the kernel $K$ given by (\ref {lblzker}).
Combining (\ref {lblzdrft}) and (\ref {lblzcll}), letting $h\rightarrow 0,\ N\rightarrow \infty $ gives
\begin {align} \label {lblzwk}
& \int \phi (x,\kappa )\partial _s f(x,\kappa ,s)\ dx\kappa =
\\ & 
\int  \nabla _2\phi (y,r)\cdot \omega _1(y,r,s) f(y,r,s) \ dyr
+
\int  [\phi (x,r)-\phi (y,r) ]
K(x,y,r)f(y,r,s) \ dxyr
\end{align}
which is the weak form of 
\begin{align}
&	\partial _s f(x,\kappa ,s)=
\\ &
-div _\kappa [\omega _1(x,\kappa ,s) f]
+
\int 
K(x,y,\kappa )f(y,\kappa ,s) \ dy
-
\int 
K(y,x,\kappa )f(x,\kappa ,s) \ dy
\end{align}
with the kernel $K$ given by (\ref {lblzker}).
\endproof
\par

{\bf Proof of Proposition \ref {prpcnv} :}\\ 
We define the symmetrized kernel $K_{sym}$ by
$K_{sym}(x,y)=K(x,y)f_\infty (y)=\alpha (x,y)\tau (x|y)f_\infty (y)$ satisfying, according to the detail balance condition (\ref{ldbblz}), the symmetry
$K_{sym }(x,y)=K_{sym}(y,x)$.
This makes (\ref {lblzwk1}) into 
$$
\int  \phi (x )\partial _s f(x,s )\ dx =
\int  [\phi (x)-\phi (y)]K_{sym}(x,y)\frac {f(y,s)}{f_\infty  (y)} \ dxy 
\ \forall \phi \ .
$$ 
Interchanging the integration variables $x\leftrightarrow y$ in the integral on the right hand side gives
(because of the symmetry of $K_{sym}$)
$$
\int  \phi (x )\partial _s f(x,s )\ dx =
\int  [\phi (y)-\phi (x)]K_{sym}(x,y)\frac {f(x,s)}{f_\infty  (x)} \ dxy \ .
\ \forall \phi \ ,
$$ 
Summing these two equations gives
\begin {equation} \label {lblzsm}
2\int  \phi (x )\partial _s f(x,s )\ dx =
\int  [\phi (x)-\phi (y)]K_{sym}(x,y)[\frac {f(y,s)}{f_\infty  (y)}-\frac {f(x,s)}{f_\infty  (x)}\ \ dxy \ .
\ \forall \phi \ .
\end {equation} 
First, we see from (\ref {lblzsm}) that we obtain a steady state for $f(x,s)=f_\infty(x)$ since the right hand side vanishes for all test functions $\phi $.
Second, choosing the special test function $\phi =\frac {f(x,s)}{f_\infty (x)}$ 
$$
2\int  \frac {f(x,s)}{f_\infty (x)}\partial _s f(x,s )\ dx =
\partial _s\int  \frac {f(x,s)^2}{f_\infty (x)}\ dx =
-\int K_{sym}(x,y)[\frac {f(y,s)}{f_\infty  (y)}-\frac {f(x,s)}{f_\infty  (x)}]^2\ \ dxy\leq 0 \ .
$$
Thus, the convex entropy functional $H_f(s)=\int  \frac {f(x,s)^2}{f_\infty (x)}\ dx $ decays monotonically until the limit $f=f_\infty $ is reached. 
\endproof 
\par 
{\bf Proof of Proposition \ref {prpbm}}
\\ 
We separate  (\ref {lcntnrm}) into a term which reduces to the pure Browninan motion (for $\beta =0$) and a term dependent on the rejection rate $\beta $.
\begin{align}
& \int  \phi (x,\kappa )f(x,\kappa ,s+h)\ dx\kappa =A+B
\\
& \label {lAbm}
A=\int  \phi (x,\omega (x,r,s) )\frac 1 {(\sqrt h \sigma  (y,r))^d} \psi (\frac {x-y-hE(y,r)}{\sqrt h \sigma  (y,r)})f(y,r,s)\ dyrx
\\ & 
 \label {lBbm}
B=\int  [ \phi (y,\omega (y,r,s) )-\phi (x,\omega (x,r,s))]\beta (x,\omega (x,r,s),y,\omega (y,r,s))\times 
\\ & 
\frac {1}{(\sqrt h \sigma  (y,r))^d}\psi (\frac {x-y-hE(y,r)}{\sqrt h \sigma  (y,r)})f(y,r,s)\ dyrx \ .
\end{align}
This gives for the term $A$ in (\ref {lAbm}),
after the variable transformation $x\rightarrow y+hE+\sqrt h \sigma  x,\ dx\rightarrow \sigma  ^d dx$,
\begin{equation} \label {lAnm}
A=\int  \phi (y+hE(y,r)+\sqrt h \sigma  (y,r)x,\omega (y+hE+\sqrt h \sigma  x,r,s) ) \psi (x)f(y,r,s)\ dyrx
\end{equation}  
The term $A$ should  (for $\beta =0$) give the classical Brownian motion term.
Similarly, for the term $B$ in (\ref {lBbm}) (the correction if we do not always accept) we obtain, normalizing $\tau $ and transforming in the integral
\begin{equation} \label {lBnm}
B=
[\phi (y,\omega (y,r,s) )-\phi (y+hE+\sqrt h \sigma  x ,\omega (y+hE+\sqrt h  \sigma  x,r,s))]\times 
\end{equation}  
$$
\beta (y+h E+\sqrt h \sigma  x,\omega (y+h E+\sqrt h  \sigma  x,r,s),y,\omega (y,r,s)) \psi (x)
f(y,r,s)\ dyrx,\
$$
We will still have to move the zero order term in $A$ to the left hand side and divide by the stepsize $h$ to obtain a derivative in $s$.
So, what remains is to expand  (\ref {lAnm}) and (\ref {lBnm}) for small stepsizes $h$ up to $o(h)$.
\par 
{\bf Expansion of $A$}.  We have from (\ref {lAnm})
$$
A=\int  \phi (y+hE(y,r)+\sqrt h \sigma  (y,r)x,\omega (y+hE+\sqrt h \sigma  x,r,s) ) \psi (x)f(y,r,s)\ dyrx
$$
We first expand $\phi $ in terms of the perturbation in the first variable $y$.
\begin{align}
& A=A_0+A_1+A_2+o(h)\ .
\\ & 
A_0=\int  \phi (y,\omega (y+hE+\sqrt h \sigma  x,r,s) ) \psi (x)f(y,r,s)\ dyrx \ .
\\ & A_1=
\int(hE(y,r)+\sqrt h \sigma  (y,r)x) \partial _1\phi (y,\omega (y+hE+\sqrt h \sigma  x,r,s) ) \psi (x)f(y,r,s)\ dyrx
\\ & A_2=
\int  \frac 1 2 (hE(y,r)+\sqrt h \sigma  (y,r)x)^2\partial _1^2\phi (y,\omega (y+hE+\sqrt h \sigma  x,r,s) ) \psi (x)f(y,r,s)\ dyrx
\end{align}
First we write the proposed moments $\omega (x,\kappa ,s)$ in (\ref {lBMom}) as a $O(h)$ perturbation around the current moments $\kappa $: 
\begin{equation} \label {lombm}
\omega (x,\kappa ,s)=\kappa +h\omega _1 (x,\kappa ,s),\
\omega _1(x,\kappa ,s)=\frac {\Delta \kappa (x)-\kappa }{s+h} \ .
\end {equation}
expanding the growth $\omega $ in the moments using $\omega (y,r,s)=r+h\omega _1(y,r,s)$
\begin{align}
	&
A_0=\int  \phi (y,r+h\omega _1(y+hE+\sqrt h \sigma  x,r,s) ) \psi (x)f(y,r,s)\ dyrx
\\ & 
=\int  \phi (y,r)  \psi (x)f(y,r,s)\ dyrx+
\int  (h\omega _1(y+hE+\sqrt h \sigma  x,r,s)\cdot  \nabla _r\phi (y,r) ) \psi (x)f(y,r,s)\ dyrx+O(h^2)
\\ & 
=\int  \phi (y,r)  \psi (x)f(y,r,s)\ dyrx+
\int  h\omega _1(y,r,s)\cdot  \nabla _r\phi (y,r)  \psi (x)f(y,r,s)\ dyrx+o(h)
\end{align} 
integrating against the normalized distribution $\psi $ gives
$$
A_0=\int  \phi (y,r)  f(y,r,s)\ dyr+
\int  h\omega _1(y,r,s)\cdot  \nabla _r\phi (y,r)  f(y,r,s)\ dyr+o(h) \ .
$$
We proceed in the same way with the terms $A_1$ and $A_2$:
\begin{align}
	&
A_1=
\int(hE(y,r)+\sqrt h \sigma  (y,r)x) \cdot  \nabla _y\phi (y,r+h\omega _1(y+hE+\sqrt h \sigma  x,r,s) ) \psi (x)f(y,r,s)\ dyrx
\\ & 
=
\int(hE(y,r)+\sqrt h \sigma  (y,r)x)\cdot  \nabla _y\phi (y,r) ) \psi (x)f(y,r,s)\ dyrx+O(h^{3/2})
\\ & 
=
\int  hE(y,r)\cdot  \nabla _y\phi (y,r) f(y,r,s)\ dyr+O(h^{3/2}) \ .
\end{align}
For $A_2$ we obtain
\begin{align}
& A_2=
\int  \frac 1 2 (hE(y,r)+\sqrt h \sigma  (y,r)x)^2\partial _1^2\phi (y,r+h\omega _1(y+hE+\sqrt h \sigma  x,r,s) ) \psi (x)f(y,r,s)\ dyrx
\\ &=
\int  \frac 1 2 (hE(y,r)+\sqrt h \sigma  (y,r)x)^2\partial _1^2\phi (y,r) \psi (x)f(y,r,s)\ dyrx+O(h^2)
\\ & =
\frac 1 2 h\int   \sigma  (y,r)^2( x)^2\partial _1^2\phi (y,r) \psi (x)f(y,r,s)\ dyrx+O(h^{3/2})
\\ & =
\frac 1 2 h\int   \sigma  (y,r)^2x_jx_k\partial _{y_jy_k}\phi (y,r) \psi (x)f(y,r,s)\ dyrx+O(h^{3/2})
\\ & =
\frac 1 2 h\int   \sigma  (y,r)^2f(y,r,s)div _y\nabla _y\phi (y,r) \ dyr+O(h^{3/2})
\end{align}
so, we have for the pure Brownian motion term (with a zero rejection rate $\beta =0$)
$A=A_0 +A_1+A_2+o(h)$ with
\begin{align} \label {lAexp}
& A =
\int  \phi (y,r)  f(y,r,s)\ dyr+
\int  h\omega _1(y,r,s)\cdot  \nabla _r\phi (y,r)  f(y,r,s)\ dyr
\\ 
&
+\int  hE(y,r)\cdot  \nabla _y\phi (y,r) f(y,r,s)\ dyr
+\frac 1 2 h\int   \sigma  (y,r)^2f(y,r,s)div _y\nabla _y\phi (y,r) \ dyr+o(h) \ .
\end{align}

{\bf Expansion of $B$:}  We define the difference in the argument of the test function in (\ref {lBnm}) as
$$
\Delta \phi (x,y,r,s)=[\phi (y,\omega (y,r,s) )-\phi (y+hE+\sqrt h \sigma  x ,\omega (y+hE+\sqrt h  \sigma  x,r,s))]
$$
with $E=E(y,r)$ and $\sigma  =\sigma  (y,r)$
and, from (\ref {lBnm}) we write
$$
B=\int  
\Delta \phi  (x,y,r,s)
\beta (y+h E+\sqrt h \sigma  x,\omega (y+h E+\sqrt h  \sigma  x,r,s),y,\omega (y,r,s)) \psi (x)
f(y,r,s)\ dyrx \ .
$$
Again, we use $\omega (y,r,s)=r+h\omega _1(y,r,s)$.
Expanding the shift in $y$ in $\Delta \phi $ gives
\begin{align}
& \Delta \phi (x,y,r,s)=
[\phi (y,\omega (y,r,s) )-\phi (y+hE+\sqrt h \sigma  x ,\omega (y+hE+\sqrt h  \sigma  x,r,s))]
\\ & 
=T_0+T_1+T_2+O(h^{3/2}) \mbox{ for }
\\ & 
T_0=
\phi (y,\omega (y,r,s) )
-\phi (y,\omega (y+hE+\sqrt h  \sigma  x,r,s))
\\ & 
T_1=
-(hE+\sqrt h \sigma  x)\cdot \nabla _1\phi (y ,\omega (y+hE+\sqrt h  \sigma  x,r,s))
\\ & 
T_2=\frac 1 2 (hE+\sqrt h \sigma  x)^2\cdot \partial _1^2\phi (y ,\omega (y+hE+\sqrt h  \sigma  x,r,s))
\end{align}
Again, expanding $\omega $ using $\omega (y,r,s)=r+h\omega _1(y,r,s)$ gives for $T_0,T_1,T_2$
\begin{align}
T_0 &=
\phi (y,r+h\omega _1(y,r,s) )
-\phi (y,r+h\omega _1(y+hE+\sqrt h  \sigma  x,r,s))
\\ &=
\phi (y,r )
-\phi (y,r)
\\ &=
+h\omega _1(y,r,s)\partial _2\phi (y,r)
-h\omega _1(y+hE+\sqrt h  \sigma  x)\partial _2\phi (y,r)
=O(h^{3/2}) \ ,
\end{align}
and for $T_1$ and $T_2,$ respectively: 
\begin{align}
	T_1 &=
-(hE+\sqrt h \sigma  x)\cdot \nabla _1\phi (y ,r+h\omega _1(y+hE+\sqrt h  \sigma  x,r,s))=
-(hE+\sqrt h \sigma  x)\partial _1\phi (y ,r) +O(h^{3/2}) \ ,
\\  
T_2 &=\frac 1 2 (hE+\sqrt h \sigma  x)^2\cdot \partial _1^2\phi (y ,r+h\omega _1(y+hE+\sqrt h  \sigma  x,r,s))=
\frac 1 2 (hE+\sqrt h \sigma  x)^2\cdot \partial _1^2\phi (y ,r)+O(h^2) \ .
\end{align}
So, altogether, we have
$$
\Delta \phi (x,y,r,s)=
-(hE+\sqrt h \sigma  x)\cdot \nabla _1\phi (y ,r)
+
\frac 1 2 (hE+\sqrt h \sigma  x)^2\cdot \partial _1^2\phi (y ,r)+O(h^{3/2})
$$
or
\begin {equation}\label {ldphiexp}
\Delta \phi (x,y,r,s)=
-(hE+\sqrt h \sigma  x)\cdot \nabla _1\phi (y ,r)
+
\frac 1 2  h \sigma  ^2x^2\cdot \partial _1^2\phi (y ,r)+O(h^{3/2})
\end {equation}
\par
{\bf Expansion of $\beta $.}  We note, that $\Delta \phi =O(\sqrt h)$.
So, we have to expand $\beta $ only up to $o(\sqrt h)$ to obtain an expression for the product $\Delta \phi \beta $ up to order $o(h)$. 
We first expand again in the state variable:
\begin{align}
& \beta (y+h E+\sqrt h \sigma  x,\omega (y+h E+\sqrt h  \sigma  x,r,s),y,\omega (y,r,s))
=\beta _0+\beta _1+o(\sqrt h), \\ 
& \beta _0=\beta (y,\omega (y+h E+\sqrt h  \sigma  x,r,s),y,\omega (y,r,s))
\\ 
& \beta _1=(h E+\sqrt h \sigma  x)\partial _1\beta (y,\omega (y+h E+\sqrt h  \sigma  x,r,s),y,\omega (y,r,s))
\end{align}
Expansion in the moment corrections $\omega $ gives,
using $\omega (y,r,s)=r+h\omega _1(y,r,s)$
\begin{align}
	&
\beta _0=\beta (y,r+h\omega _1(y+h E+\sqrt h  \sigma  x,r,s),y,r+h \omega _1(y,r,s))=
\beta (y,r,y,r)+o(\sqrt h)
\\ 
& \beta _1=(h E+\sqrt h \sigma  x)\cdot \nabla _1\beta (y,r+h\omega _1(y+h E+\sqrt h  \sigma  x,r,s),y,r+h\omega _1(y,r,s))=
\\ 
& 
(h E+\sqrt h \sigma  x)\cdot \nabla  _1\beta (y,r,y,r)+o(\sqrt h)
\end{align}
Therefore, we have
$$
\beta (y+h E+\sqrt h \sigma  x,\omega (y+h E+\sqrt h  \sigma  x,r,s),y,\omega (y,r,s))=
$$
\begin {equation} \label {lbetaexp}
\beta (y,r,y,r)+(h E+\sqrt h \sigma  x)\partial _1\beta (y,r,y,r)+o(\sqrt h)
\end {equation}
multiplying $\Delta \phi $ with $\beta $ from (\ref {ldphiexp}) and (\ref {lbetaexp}) gives
\begin{align}
& \Delta \phi \beta =
[-(hE+\sqrt h \sigma  x)\cdot \nabla _1\phi (y ,r)
+
\frac 1 2  h \sigma  ^2x^2\cdot \partial _1^2\phi (y ,r)+O(h^{3/2})]\times 
\\ & 
[\beta (y,r,y,r)+(h E+\sqrt h \sigma  x)\cdot \nabla _1\beta (y,r,y,r)+o(\sqrt h)]
\\ & 
\Delta \phi \beta =
-(hE+\sqrt h \sigma  x)\cdot \nabla _1\phi (y ,r)\beta (y,r,y,r)
+\frac 1 2  h \sigma  ^2x^2\cdot \partial _1^2\phi (y ,r)\beta (y,r,y,r)
\\ & 
-(hE+\sqrt h \sigma  x)\cdot \nabla _1\phi (y ,r)(h E+\sqrt h \sigma  x)\partial _1\beta (y,r,y,r)
\\ & 
+\frac 1 2  h \sigma  ^2x^2\partial _1^2\phi (y ,r)(h E+\sqrt h \sigma  x)\partial _1\beta (y,r,y,r)
+o(h)
\end{align}
neglecting the $o(h)$ terms gives
\begin{align}
& \Delta \phi \beta =
-(hE+\sqrt h \sigma  x)\cdot \nabla _1\phi (y ,r)\beta (y,r,y,r)
+\frac 1 2  h \sigma  ^2x^2\cdot \partial _1^2\phi (y ,r)\beta (y,r,y,r)
\\ & 
-h[ \sigma  x\partial _1\phi (y ,r)][\sigma  x\partial _1\beta (y,r,y,r)]
+o(h)
\end{align}
integrating against the normalized distribution $\psi $,
using $\int  (1,x,xx^T)\psi (x)\ dx=(1,0, \; Id)$ gives
\begin{align}
& \int  \Delta \phi \beta \psi (x)\ dx=
-hE\cdot  \nabla _y\phi (y ,r)\beta (y,r,y,r)
+\frac 1 2  h \sigma  ^2\beta (y,r,y,r)div_y\nabla _y\phi (y ,r)
\\
&
-h\sigma  ^2\nabla _y\phi (y ,r)\cdot  \nabla _1\beta (y,r,y,r)
+o(h)
\end{align}
multiplying with $f(y,r)$ and integrating $dyr$ give the term $B$ in (\ref {lBnm}) up to terms of order $o(h)$.
\begin{align} \label {lBexp}
B&=
\int  [-hE\cdot  \nabla _y\phi (y ,r)\beta (y,r,y,r)
+\frac 1 2  h \sigma  ^2\beta (y,r,y,r)div_y\nabla _y\phi (y ,r) \\
& -h\sigma  ^2\nabla _y\phi (y ,r)\cdot  \nabla _1\beta (y,r,y,r)] f(y,r,s)\ dyr
+o(h)
\end{align}  
combining (\ref {lAexp}) and (\ref {lBexp}) gives
$$
\int  \phi (x,\kappa )f(x,\kappa ,s+h)\ dx\kappa =A+B +o(h)
$$
or, after changing the variables $y\rightarrow x$ and $r\rightarrow \kappa $ in the integrals  on the right hand side
\begin{align}
& \int  \phi (x,\kappa )f(x,\kappa ,s+h)\ dx\kappa =
\\ & 
\int  \phi (x,\kappa )  f(x,\kappa ,s)\ dx\kappa +
\int  h\omega _1(x,\kappa ,s)\cdot  \nabla _\kappa \phi (x,\kappa )  f(x,\kappa ,s)\ dx\kappa 
\\ & 
+\int  hE(x,\kappa )\cdot  \nabla _x\phi (x,\kappa ) f(x,\kappa ,s)\ dx\kappa 
+\frac 1 2 h\int   \sigma  (x,\kappa )^2f(x,\kappa ,s)div _x\nabla _x\phi (x,\kappa ) \ dx\kappa 
\\ & 
+
\int  [-hE\cdot  \nabla _x\phi (x ,\kappa )\beta (x,\kappa ,x,\kappa )
+\frac 1 2  h \sigma  ^2\beta (x,\kappa ,x,\kappa )div_x\nabla _x\phi (x,\kappa )
\\ & 
-h\sigma  ^2\nabla _x\phi (x ,\kappa )\cdot  \nabla _1\beta (x,\kappa ,x,\kappa )] f(x,\kappa ,s)\ dx\kappa 
+o(h)
\end{align}
moving the $O(1)$ term $ \int  \phi (x,\kappa )  f(x,\kappa ,s)\ dx\kappa  $ to the left of the equality sign  and dividing by the stepsize $h$, and letting $h\rightarrow 0$ gives for all test functions $\phi:$ 
\begin{align}
	& 
\int  \phi (x,\kappa )\partial _sf(x,\kappa ,s)\ dx\kappa =
\\ & 
\int  \omega _1(x,\kappa ,s)\cdot  \nabla _\kappa \phi (x,\kappa )  f(x,\kappa ,s)\ dx\kappa 
\\ & 
+\int  E(x,\kappa )\cdot  \nabla _x\phi (x,\kappa ) f(x,\kappa ,s)\ dx\kappa 
+\frac 1 2 \int   \sigma  (x,\kappa )^2f(x,\kappa ,s)div _x\nabla _x\phi (x,\kappa ) \ dx\kappa 
\\ & 
+
\int  [-E\cdot  \nabla _x\phi (x ,\kappa )\beta (x,\kappa ,x,\kappa )
+\frac 1 2   \sigma  ^2\beta (x,\kappa ,x,\kappa )div_x\nabla _x\phi (x,\kappa )
\\ & 
-\sigma  ^2\nabla _x\phi (x ,\kappa )\cdot  \nabla _1\beta (x,\kappa ,x,\kappa )] f(x,\kappa ,s)\ dx\kappa 
+o(1)
\end{align}
After integrating by parts this is the weak form of
\begin{align}
& \partial _sf(x,\kappa ,s)=
-div _\kappa [\omega _1(x,\kappa ,s)f(x,\kappa ,s)]
\\ & 
-div_x[ E(x,\kappa ) f(x,\kappa ,s)]
+\frac 1 2 div _x \nabla _x [ \sigma  (x,\kappa )^2f(x,\kappa ,s)]
\\ & 
+div _x[\sigma  ^2\nabla _1\beta (x,\kappa ,x,\kappa ) f(x,\kappa ,s)] \ .
\end{align}
We consolidate
\begin{align} \label {ldsf}
& \partial _sf(x,\kappa ,s)=
-div _\kappa [\omega _1(x,\kappa ,s)f(x,\kappa ,s)]
\\ & 
+div_x[ E(x,\kappa ) (\beta (x,\kappa ,x,\kappa )-1)f(x,\kappa ,s)+\sigma  ^2\nabla _1\beta (x,\kappa ,x,\kappa ) f(x,\kappa ,s)]
\\ & 
+\frac 1 2 div _x \nabla _x [ \sigma  (x,\kappa )^2(1+\beta (x,\kappa ,x,\kappa ))f(x,\kappa ,s)]
\end{align}
\endproof 

\subsection{Implementation details} \label{details}
We report on the  details required for the simulation of the Metropolis Monte Carlo algorithm of Section \ref{ssgendef}. Here, all parameter choices are detailed. For readability we present a  the parameters in the  order of appearance.
\par \noindent
 The total number of time steps $n=1,\dots,N$ is $N=10'000$ in all simulations.
\par\noindent \noindent
In all cases we set $T=25$,  $K_{\max}=10'000$ and  $a_1=a_3=10$ and $a_2=1.$
	\par\noindent The moments $\kappa$ defined by equation \eqref{kappa} we consider $\kappa=(\kappa_1,\sigma)$ defined by 
	\begin{align}
		\kappa^n _1 = \frac1n \sum_{j=1}^n  x_j, \; 
		\sigma^n     =  \kappa^n_2 - (\kappa_1^n)^2, \kappa_2^n = \frac1n \sum_{i=1}^n   \left( x_j  \right)^2.
	\end{align} 
	\par\noindent 
 According to remark \eqref{rss}, the Monte Carlo algorithm requires a start up phase. We realize this by producing an initial distribution obtained by considering $i=1,\dots, N$ simulations of model \eqref{v63} with parameter
	$x_i:=x^* + \xi_i$ and  $\xi_i$ being a realization of the random variable $\xi$ distributed according to equation \eqref{v63xi}. This yields the distribution 
	\begin{align}
		P_0(x) = \frac1{N_0} \sum_{i=1}^{N_0} \delta(x - x_i),
	\end{align}
	and the corresponding moments $\kappa_0.$ In our simulations we set ${N_0}=100.$
	\par\noindent 
	 The ODE system \eqref{v63} is solved on the time interval $[0,T]$ using an  explicit 3(2) Runge-Kutta of the Matlab routine $ode23.$
	\par\noindent 
	 The acceptance rate $\alpha$ is chosen independently of the update $\omega$ in all simulations. We set 
	\begin{align}
		\alpha(x_p,x_n)=\frac{ L(x_p, \nu) }{L(x_p,\nu) + L(x_n,\nu )},
	\end{align}
	for a likelihood function $L.$ The likelihood function depends on the test case. In the first example, we set at iteration $n$  
	\begin{align}\label{fixedtimeL}
		L(x,\nu) = \frac{1}{  \| v(T,x) - z_n \|^2  },
	\end{align}
	while in the second example, we set 
	\begin{align}
		L(x,\nu) = \frac{1}{ \frac1n \sum_{i=1}^n \| v(t_i,x) - z_i \|^2  },
	\end{align}
	\par\noindent 
	 The probability for the random proposal $\tau(z | x_n, \kappa_n )$ is chosen either as Gaussian distribution or using a gradient approach. In case of the  Gaussian distribution we consider a random variable $Z$ with  distribution as 
	\begin{align}\label{34}
		Z \sim \mathcal{N}(x_n, Id_{3\times 3}), 
	\end{align} 
	and obtain the proposal $x_p$ as 
	$x_p = P_{[ \frac54 x^*, \frac34 x^* ]} ( \tilde{Z} ),$ where $\tilde{Z}$ is a realisation of $Z$ and $P_{[a,b]}(x)$ denotes the projection of $x$ on the interval $[a,b].$ The later is necessary, since for $T$ large the system \eqref{v63} is known to exhibit chaotic behavior and the solution may diverge for proposals $x_p$ far from $x^*.$

	In the case of a gradient based approach we proceed as follows. With probability $p \sim \frac12$ we set 
	$x_p = P_{[ \frac54 x^*, \frac34 x^* ]} \left( x_n + \tilde{Z}  \right)$ or $x_p = P_{[ \frac54 x^*, \frac34 x^* ]} \left( x_n - \tilde{Z} \right),$ respectively.   The value of $\tilde{Z}$ is obtained as (an approximation) to the  gradient of the likelihood with respect to the parameters $\nabla_x L(x_n,\nu).$ Evaluation of the gradient requires to differentiate the solution $v(t^*,x)$ with respect to $x$ for some fixed $t^*.$ Due to equation \eqref{v63}, the variations of $v$ with respect to $x$ fulfill a  linear ODE system of dimension $3\times 3$ given by 
	\begin{align}
		\frac{d}{dt}\left(	\frac{ \partial v}{\partial x}  \right) = \begin{pmatrix}
			v_2 - v_1 + a ( \partial_a v_2(t) - \partial_a v_1(t) ) & a (\partial_b v_2 - \partial_b v_1 ) & a (\partial_c v_2 - \partial_c v_1) \\
			- v_1 - a \partial_a v_1 - \partial_a v_2 - \partial_a (v_1 v_3)  & -a \partial_b v_1 - \partial_b (v_2-v_1 v_3) &  -a \partial_c v_1  \partial_c (v_2-v_1 v_3)  \\
			\partial_a (v_1 v_3 - b v_3) - b &   \partial_b (v_1 v_3 - b v_3) - v_3 - (c+a)  &  \partial_c (v_1 v_3 - b v_3) - b \\
		\end{pmatrix}
	\end{align}
	In the previous equation we omit the dependence on $t.$ The initial conditions are $\frac{ \partial v}{\partial x}(t=0)  = 0_{3 \times 3}.$ Solving this system for each proposal is computationally prohibitive. We therefore use an explicit Euler discretization with a single time step $t^*$ to approximate it's solution. This leads to 
	\begin{align}
		\frac{ \partial v}{\partial x}(t^*,x) \approx t^* \begin{pmatrix} 	v_2 - v_1  & 0 & 0  \\
			- v_1  & 0 &  0  \\
			- b &    - (c+a)  &   - b \\
		\end{pmatrix}.
	\end{align}  
	In the case of fixed terminal time $L$ is given by \eqref{fixedtimeL}, we have $t^*=T$,  and  obtain the following  explicit form of $\tilde{Z}$ at iteration $n:$
	\begin{align}
		\tilde{Z}_i = - \frac{2}{ \| v(T,x_n) - z_n \|^4 }\left( v(T,x_n) - z_n \right) \cdot \frac{ \partial v_i }{ \partial x}(T,x_n), \;    i=1,\dots,3.
	\end{align}
	
	\par\noindent 
	 For the micro--macro decomposition we need to specify the distribution parameter $\zeta_n.$ 	
	To this end, we compute the maximal variance of the data $\nu$ up to data point $n,$  
	$	\sigma_\nu^2 = \| \mathbb{V}ar(z_1,\dots,z_n) \|_\infty, $
	the variance  $ \mathbb{V}ar((f_{macro})_{n-1})$ of the prior   macroscopic approximation $(f_{macro})_{n-1}$, 
	and the variance of the current parameter distribution 
	$\sigma_n^2 =  \mathbb{V}ar(P_n).$
	Depending on the  ratio $r$
	\begin{align}\label{rate-r}
		r = \frac{\sigma_\nu^2 - \mathbb{V}ar((f_{macro})_{n-1})}{ \sigma_n^2}  
	\end{align}
	we update $\zeta_n$ as 
	\begin{align}
		\zeta_n = \begin{pmatrix}
			P_{[0,1]}(\sqrt{r})   & r > 0 \\
			\zeta_n{n-1} & r<0. 
		\end{pmatrix}
	\end{align}
	
	\par\noindent 
		In case of the micro--macro decomposition, we  decide according to $\zeta_n$ on the probability density $\tau(z | x_n, \kappa_n )$, i.e., 
	if the realisation of a uniform random variable $Z \sim \mathcal{U}(0,1)$ at iteration $n$ is less than $\zeta_{n-1}$ we sample from the microscopic distribution. 	 In this case,  $\tau$ is independent of $\kappa_n$ and the sampling is described above in equation \eqref{34}. 
	For the  macroscopic case, the probability density $\tau$ is independent of $x_n$ and solely depends on the moments $\kappa_n$ of the macroscopic distribution $g_n.$ The proposal $x_p$ is then given by $x_p =  P_{[ \frac54 x^*, \frac34 x^* ]} \left(  \tilde{Z}  \right)$ where $\tilde{Z}$ is a realization of the random variable $Z$ 
	\begin{align}\label{35}
		Z \sim   \mathcal{N}( \mathbb{E}(g_{n-1}),  \mathbb{V}ar((f_{macro})_{n-1}) ), 
	\end{align} 
	
\end{document}